\documentclass[11pt]{amsart}

\usepackage{amssymb, amsmath, amscd, amsthm, graphicx, psfrag}
\usepackage[matrix,arrow,curve]{xy}
\xyoption{dvips}              

\newcommand{\T}{\Psi}

\newcommand{\R}{{\mathbb R}}
\newcommand{\Z}{{\mathbb Z}}
\newcommand{\Q}{{\mathbb Q}}

\newcommand{\calD}{{\mathcal{D}}}
\newcommand{\calP}{{\mathcal{P}}}

\newcommand{\nnn}{\underline{n}}

\newcommand\rth{\refstepcounter{equation}}
\newcommand\numb{\rth{\rm \theequation}}
\numberwithin{equation}{section}

\DeclareMathOperator*{\holim}{holim}

\newcommand{\BR}{{\mathbb R}}

\newcommand{\Top}{\mathrm{Top}}

\newcommand{\calC}{{\mathcal{C}}}

\newcommand{\calF}{{\mathcal{F}}}

\theoremstyle{plain}
\newtheorem{thm}{Theorem}[section]
\newtheorem{theorem}{Theorem}

\newtheorem{prop}[thm]{Proposition}
\newtheorem{proposition}[thm]{Proposition}
\newtheorem{lemma}[thm]{Lemma}
\newtheorem{cor}[thm]{Corollary}

\newtheorem{theorem-assertion}[thm]{Theorem}
\newtheorem*{theorem*}{Theorem}

\theoremstyle{definition}
\newtheorem{definition}[thm]{Definition}
\newtheorem{remark}[thm]{Remark}

\theoremstyle{remark}

\begin{document}


\title{Associahedron, Cyclohedron, and Permutohedron \\
as compactifications of configuration spaces}


\author{Pascal Lambrechts}
\address{Institut Math\'{e}matique, 2 Chemin du Cyclotron, B-1348 Louvain-la-Neuve, Belgium}
\email{lambrechts@math.ucl.ac.be}
\urladdr{http://milnor.math.ucl.ac.be/plwiki}
\author{Victor Turchin}
\address{Kansas State University, USA.}
\email{turchin@ksu.edu} \urladdr{http://www.math.ksu.edu/\~{}turchin/}
\author{Ismar Voli\'c}
\address{Department of Mathematics, Wellesley College,
Wellesley, MA}
\email{ivolic@wellesley.edu}
\urladdr{http://palmer.wellesley.edu/\~{}ivolic}

\thanks{The first author is a chercheur qualifi\'e au F.N.R.S.  The third author was supported in part by the National Science Foundation grant DMS 0504390.  }


\begin{abstract}
As in the case of the associahedron and cyclohedron, the permutohedron can also be defined as an appropriate compactification of a configuration space of points on an interval or on a circle. The construction of the compactification endows the permutohedron with a projection to the cyclohedron, and the cyclohedron with a projection to the associahedron.
We  show that the preimages of any point via these
projections might not be homeomorphic to (a cell decomposition of) a
disk, but are still contractible.  We briefly explain an application
of this result to the study of knot spaces from the point of view of
the Goodwillie-Weiss manifold calculus.
\end{abstract}

\maketitle

\sloppy

\markright{ASSOCIAHEDRON, CYCLOHEDRON, AND PERMOTOHEDRON AS COMPACTIFICATIONS}

\section{Introduction}\label{s1}

\begin{figure}[h!]
\includegraphics[width=15cm]{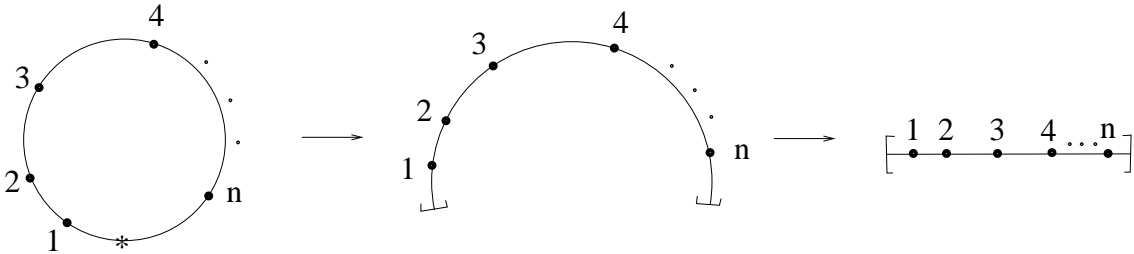}
\caption{A homeomorphism between $Conf_*(n,S^1)$ and $Conf(n,[0,1])$.}\label{homeofig}
\end{figure}

The configuration space $Conf(n,[0,1])$ of $n$ distinct points $0<t_1<t_2<\ldots<t_n<1$
in the interior of the segment $[0,1]$ is clearly homeomorphic to the configuration
space $Conf_*(n,S^1)$ of $n+1$ distinct points on the circle $S^1\simeq [0,1]/0{\sim}1$, one of which is the fixed point $*=0{\sim}1$.  This homeomorphism is pictured in Figure \ref{homeofig}.

The Axelrod-Singer compactification~\cite{AS} of $Conf(n,[0,1])$  is the $n$-dimensional associahedron $Assoc_n$, also called the Stasheff polytope~\cite{Stasheff-HAH,Stasheff-HHPV}. The Axelrod-Singer compactification of $Conf_*(n,S^1)$ is the $n$-dimensional cyclohedron $Cycl_n$, also called the Bott-Taubes polytope~\cite{BT}.

The homeomorphism $Conf_*(n,S^1)\to Conf(n,[0,1])$ induces a natural projection of compactifications
$$
\pi_n\colon Cycl_n\to Assoc_n.
\eqno(\numb)\label{eq11}
$$
The compactification $Cycl_n$ contains more information than $Assoc_n$ since
one can compare how fast configuration points approach $*=0{\sim}1$ from the left and from the right:

\vspace{0.2cm}

\begin{center}
\includegraphics[width=4cm]{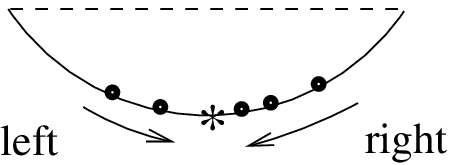}
\end{center}

\vspace{0.2cm}

The projection~\eqref{eq11} essentially forgets this information.

The poset of faces of $Assoc_n$ is the poset of planar trees (see Section~\ref{ss21}) which we denote by $\Psi([n+1])$.\footnote{Here $[n+1]$ stands for the ordered set $\{0,1,\ldots,n+1\}$
of leaves.}  We denote the poset of faces of $Cycl_n$ by $\Phi(\underline{n})$.\footnote{Here $\underline{n}$ stands for the cyclically ordered pointed set $\{0=*,1,2,\ldots,n\}$.} The elements of $\Phi(\underline{n})$  are certain planar
trees that we call {\it fans} (see Section \ref{ss22}).

To any face of $Cycl_n$ one can assign a face of $Assoc_n$ which is its image via $\pi_n$. This correspondence defines a functor $\Pi_n\colon \Phi(\underline{n})\to
\Psi([n+1])$.
Passing to a map of realizations $|\Phi(\underline{n})|\to
|\Psi([n+1])|$ gives the projection $|\Pi_n|\colon Cycl_n\to Assoc_n$.
However, the projections $\pi_n$ and $|\Pi_n|$ are different because $|\Pi_n|$ is not a homeomorphism of interiors (see Remark~\ref{r110}).

We prove the following (see Section~\ref{s3}):

\vspace{0.2cm}

\noindent {\bf Theorem~\ref{t1}.}
{\it The preimage of any point of $|\Psi([n+1])|=Assoc_n$ under $|\Pi_n|$ is contractible.}

\vspace{0.2cm}

%
%
%
%
%

We also consider the initial projection~\eqref{eq11} and describe the geometry of preimages under~$\pi_n$ (see Section~\ref{s5}-\ref{s6}). In particular we prove:

\vspace{0.2cm}

\noindent {\bf Theorem~\ref{t2}.} {\it The preimage of any point of $Assoc_n$ under $\pi_n$ is contractible.}

\vspace{0.2cm}

\begin{figure}[h]
\includegraphics[width=3cm]{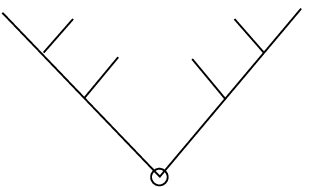}
\caption{}
\label{badtreefigure}
\end{figure}

\begin{figure}[h]
\includegraphics[width=5cm]{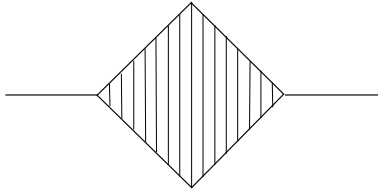}
\caption{Preimage of the tree from Figure \ref{badtreefigure}.}
\label{preimagedisk}
\end{figure}

Theorems~\ref{t1} and~\ref{t2} are not surprising, but what is interesting is that the preimages might not be homeomorphic to a disk. For example, for the vertex of $Assoc_4$ encoded by the binary tree from Figure \ref{badtreefigure} (which is the limit of $(t_1,t_2,t_3,t_4)=(\varepsilon^2,\varepsilon,1-\varepsilon, 1-\varepsilon^2)$ when $\epsilon\to +0$), the preimage under both $\pi_n$ and $|\Pi_n|$ is a square (2-disk) with two
segments attached \footnote{This preimage is the realization of the poset $X_{2,2}$; see Lemma~\ref{l32} and Proposition~\ref{p61}~(iii).} as in Figure~\ref{preimagedisk}.

\vspace{0.2cm}

\vspace{0.2cm}

We were not able to find the projection $\Pi_n$ in the literature, even though it can be easily derived from the work in ~\cite{Tonks-RAP} (and it was actually known to A.~Tonks). For other interesting relations between associahedra, cyclohedra, and other polytopes see~\cite{HohLange-RAC,Markl,Postnikov-PAB,Devadoss-Cycl}.



In Section~\ref{s8}, we define a new {\it leveled compactification} $C_n[[M]]$ of a configuration space $C_n(M)$ of $n$ points in a smooth manifold $M$. The difference from the usual Axelrod-Singer compactification $C_n[M]$ is that  the new one takes into account the ratios of diameters of distant infinitesimal conglomerations of points. This compactification can also be obtained by the construction given by Gaiffi in~\cite{Gaiffi} (see Section~\ref{s8}). One has a natural projection $C_n[[M]]\longrightarrow C_n[M]$. The $n$-dimensional permutohedron can be defined as a leveled compactification of $Conf(n,[0,1])$ or of $Conf_*(n,S^1)$. The projection from the leveled compactification to the usual one produces maps $Perm_n\stackrel{\pi''_n}{\longrightarrow}Assoc_n$, $Perm_n\stackrel{\pi'_n}{\longrightarrow}Cycl_n$. According to A.~Tonks~\cite{Tonks-RAP}, the face poset of $Perm_n$ is the poset of leveled trees $\Psi^{level}([n+1])$. This poset is naturally isomorphic to the poset of leveled fans $\Phi^{level}(\nnn)$ (see Lemma~\ref{l71}). Forgetting levels defines  projections $\Phi^{level}(\nnn)\stackrel {\Pi'_n}\longrightarrow \Phi(\nnn)$, and $\Psi^{level}([n+1])\stackrel{\Pi''_n}\longrightarrow\Psi([n+1])$. (The poset maps $\Pi'_n$, $\Pi''_n$ describe which face of $Perm_n$ is sent to which face of $Cycl_n$ or $Assoc_n$ via $\pi'_n$ or~$\pi''_n$ respectively.) Passing to the realizations, one gets the projections
$$
Perm_n\stackrel{|\Pi'_n|}\longrightarrow Cycl_n, \hspace{1cm} Perm_n\stackrel{|\Pi''_n|}\longrightarrow Assoc_n.
\eqno(\numb)\label{eq01}
$$
Similarly to Theorem~\ref{t1}, we prove the following (see Section~\ref{s7}):

\vspace{0.2cm}

\noindent {\bf Theorem~\ref{t3}.}
{\it The preimage of any point of $Cycl_n$ (respectively $Assoc_n$) under $|\Pi'_n|$ (respectively $|\Pi''_n|$) is contractible.}

\vspace{0.2cm}

Theorem~\ref{a83} describes the cellular decomposition of any point under the projection $C_n[[M]]\longrightarrow C_n[M]$. One of the consequences is that these preimages are always contractible, which implies the following (see Section~\ref{s8}):

\vspace{0.2cm}

\noindent {\bf Theorem~\ref{a4}.}
{\it The preimage of any point of $Cycl_n$ (respectively $Assoc_n$) under $\pi'_n$ (respectively $\pi''_n$) is contractible.}

\vspace{0.2cm}

An immediate corollary of Theorems~\ref{t1} and~\ref{t3} is the following (see Section~\ref{s_application}):

\vspace{0.2cm}

\noindent {\bf Theorem~\ref{t5}.}
{\it The functors $\Pi_n$, $\Pi'_n$, $\Pi''_n$ are left cofinal.}

\vspace{0.2cm}

The above result and its application to manifold calculus was our main motivation for this paper. Theorem~\ref{t5} implies equivalence of different models for the Goodwillie-Weiss embedding tower (see Section~\ref{s_application}).

\section{Categories of faces}\label{s2}

\subsection{Category of trees}\label{ss21}
In this section we define a category $\Psi([n])$ of trees which keeps track of  the faces of $Assoc_{n-1}$.

\begin{definition}\label{def-tree}
 A {\it $\T$-tree} is an isotopy class of rooted trees embedded in the upper
half-plane with the root of valence $\geq 2$ at the origin. The valence of any internal vertex (i.e. a vertex that is not a leaf) except the root is at least 3.
\end{definition}

We orient each edge of a $\T$-tree from the vertex closer to
the root to that which is farther from the root. Each vertex
(except the root) has exactly one incoming edge and a linearly
ordered (clockwise) set of outgoing edges. The root has only
outgoing edges which are linearly ordered (clockwise).

The set of leaves has a natural (clockwise) linear order. More precisely,  let $v_{1}$ and $v_{2}$ be two leaves.
Consider two paths -- one from the root to $v_{1}$, and
another from the root to $v_{2}$. Suppose
$e_{1}$ and $e_{2}$ are the first edges that are different in these paths. These edges
 are outgoing from some vertex and we say
$v_{1}<v_{2}$ if and only if $e_{1}<e_{2}$. In particular  we
can thus speak of the minimal and maximal leaf.

\begin{definition}\label{d12}
A {\it left-most (resp. right-most) node} of a $\T$-tree is any vertex lying on the path from the
root to the minimal (resp. maximal) leaf. (Neither the root, nor the extremal leaves are considered to be left-most or right-most.)
\end{definition}

\begin{definition}\label{d13}
Define $\T([n])$ as the category whose objects are $\T$-trees with $n+1$ leaves labeled by the ordered set $[n]=\{0, 1, 2, ..., n\}$.
There is a (unique) morphism in $\T([n])$ from $T$ to $T'$ if $T'$ is obtained from $T$ by a contraction along some set
of non-leaf edges.
\end{definition}

We will think of $\T([n])$ as a poset by saying $T\geq T'$ in the above situation.

\begin{remark}\label{r14}
{\rm The realization of the category $\T([n])$ is isomorphic as a simplicial complex to the barycentric subdivision of the $(n-1)$-dimensional Stasheff associahedron $Assoc_{n-1}$.}
\end{remark}

Categories $\T([1])$, $\T([2])$, $\T([3])$ are pictured on the right side of Figure~\ref{fig1}. The root is
designated by a little circle.

\subsection{Category of fans}\label{ss22}
In this section we define a category $\Phi(\underline{n})$ of fans which keeps track of the faces of $Cycl_n$.

\begin{definition}\label{d15}
A {\it fan} (or a {\it $\Phi$-tree}) is an isotopy class of planar rooted trees with one marked leaf, called distinguished leaf. The root can be of any valence $\geq 1$.
\end{definition}

Notice that fans are not supposed to lie in an upper half-plane contrary to the $\Psi$-trees. The leaves of a fan have a natural cyclic order. The root even if it has valence one does not count for a leaf.

\begin{definition}\label{d17}
Define $\Phi(\nnn)$ to be the category whose objects are fans with $n+1$ leaves, labeled by the cyclically ordered pointed set $\nnn=\{0=*,1, 2, ..., n\}$, where $*$ labels the distinguished leaf of a fan. There is a (unique) morphism in $\Phi(\nnn)$ from $\widehat T$ to $\widehat T'$ if
$\widehat T'$ is obtained from $\widehat T$ by a contraction along some set of non-leaf edges.
\end{definition}

We will think of $\Phi(\nnn)$ as a poset by saying $\widehat T \geq\widehat T'$ in the above situation.

\begin{remark}
{\rm The realization of the category $\Phi(\nnn)$ is isomorphic as a simplicial complex to the barycentric subdivision of the $n$-dimensional
cyclohedron $Cycl_n$. Indeed, in the same way as  $\Psi([n+1])$ describes linear parenthesizations on $[n+1]$, the poset $\Phi(\nnn)$ describes cyclic parenthesizations on the set $\nnn$. This polytope was introduced by R.~Bott and C.~Taubes in~\cite{BT}.}
\end{remark}

Categories $\Phi(\underline{0})$,
$\Phi(\underline{1})$,
$\Phi(\underline{2})$ are pictured on the left side of Figure~\ref{fig1}. The root is
designated by a little circle; the distinguished leaf is designated by a black point.

\subsection{Functor $\Pi_n$}\label{ss23}

We now define a functor $\Pi_{n}\colon\Phi(\nnn)\to\T([n+1])$ between the categories of fans and trees.

\begin{definition}\label{d19}
(i)  The (only) path from the root to the distinguished leaf will be called the {\it  trunk} of the fan.

(ii)  Let $v$ be a vertex on the trunk of a fan. Let $e_{1}$ (resp. $e_{2}$) be the edge of the trunk which is
adjacent to $v$ and whose other vertex is closer to the distinguished leaf (resp. root). The edges adjacent to $v$ and lying between
$e_{1}$ and $e_{2}$ (resp. $e_{2}$ and $e_{1}$) with respect to the natural clockwise cyclic order, will be called {\it
left-going branches} (resp. {\it right-going branches}).
\end{definition}

\begin{figure}[h!]
\psfrag{Th}[0][0][1][0]{$\widehat T$}
\psfrag{T}[0][0][1][0]{$T$}
\includegraphics[width=15cm]{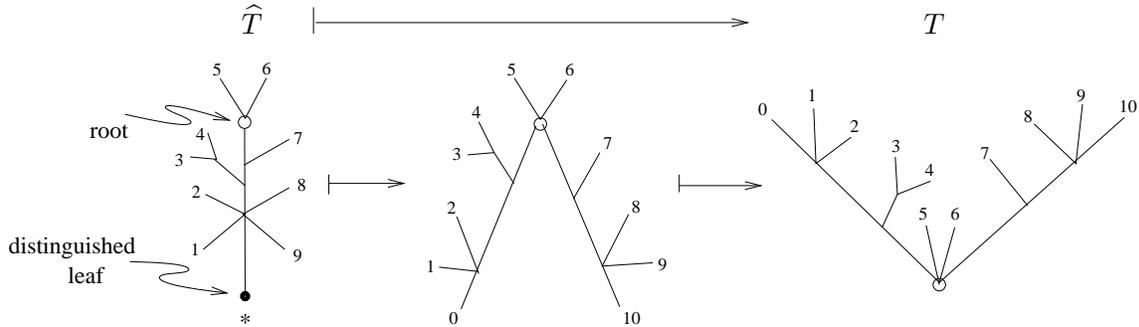}
\caption{A $\T$-tree obtained from a fan.}\label{fig0}
\end{figure}

Now let $\widehat T\in\Phi(\nnn)$ be a fan and cut $\BR^2$ along the path which is the union of the trunk of $\widehat T$ and the ray emanating downward from the distinguished leaf to infinity (so this
ray does not cross $\widehat T$). The space obtained from $\BR^2$ by this surgery is homeomorphic to the upper
half-plane. After this operation, the fan $\widehat T$ becomes a $\T$-tree $T$ with $n+2$ leaves (see Figure~\ref{fig0}).

Note that a node of $\widehat T$ along the trunk can produce either one or two vertices in $T$.
Such a node produces a left-most (resp. right-most) vertex in $T$ if and only if it has left-going (right-going) branches.
The distinguished leaf always produces two leaves, the minimal and the maximal one of the tree.


The following is immediate from the definition.

\begin{lemma}\label{l110}
The above correspondence defines a functor $\Pi_{n}\colon\Phi(\nnn)\to\T([n+1]).$
\end{lemma}

\begin{figure}
\psfrag{pi0}[0][0][1][0]{$\Pi_0$}
\psfrag{pi1}[0][0][1][0]{$\Pi_1$}
\psfrag{pi2}[0][0][1][0]{$\Pi_2$}
\psfrag{phi0}[0][0][1][0]{}
\psfrag{psi0}[0][0][1][0]{}
\psfrag{phi1}[0][0][1][0]{}
\psfrag{psi1}[0][0][1][0]{}
\psfrag{phi2}[0][0][1][0]{}
\psfrag{psi2}[0][0][1][0]{}
\psfrag{=}[0][0][1][0]{}
\includegraphics[width=16cm]{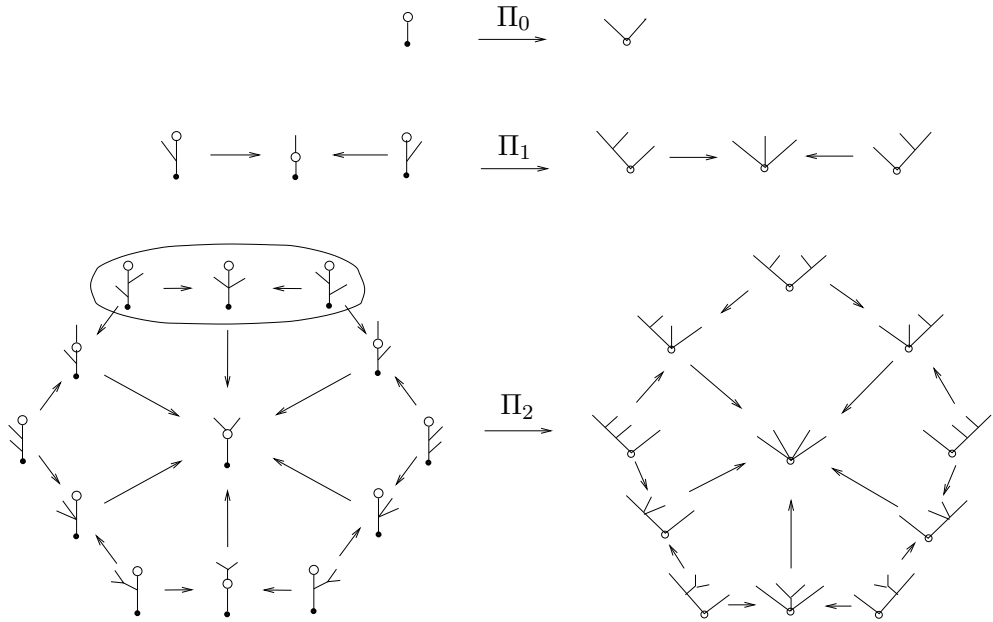}
\caption{Projections $\Pi_{0}$, $\Pi_{1}$, and $\Pi_{2}$. The circled edge of the hexagon gets mapped to the top vertex of the pentagon.}\label{fig1}
\end{figure}

\begin{remark}\label{r110}
{\rm Passing to realizations, $\Pi_n$ defines a projection $|\Pi_n|\colon|\Phi(\nnn)|\to|\T([n+1])|$. Notice however that $|\Pi_n|$ is different
from the projection $\pi_n\colon Cycl_n\to Assoc_n$ mentioned in Introduction. $|\Pi_n|$ is not a homeomorphism of interiors of $Cycl_n$ and $Assoc_n$ starting from $n\geq 2$. For example $|\Pi_2|$ in Figure~\ref{fig1} maps the 2-simplex
$\langle \includegraphics[width=.3cm]{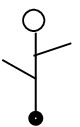}, \includegraphics[width=.3cm]{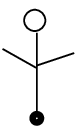}, \includegraphics[width=.3cm]{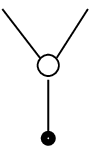} \rangle$ on the 1-simplex
$\langle \includegraphics[width=.6cm]{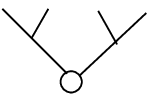} ,\includegraphics[width=.6cm]{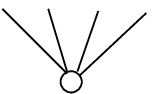} \rangle$.}
\end{remark}

\begin{theorem}\label{t1}
The preimage of any point of $|\Psi([n+1])|=Assoc_n$ under $|\Pi_n|$ is contractible.
\end{theorem}

We prove this result in the next section.

\section{Proof of Theorem~\ref{t1}}\label{s3}

To prove Theorem~\ref{t1}, we need to define certain posets which are essential for understanding the geometry of the functor
$\Pi_{n}\colon\Phi(\nnn)\to\T([n+1])$ and then prove that they are contractible.  Theorem~\ref{t1} will in the end follow from Propositions~\ref{p36} and~\ref{p_contractible}.

\begin{definition}\label{d21}
Define $X_{\ell,r}$, $l,r\geq 0$, to be the poset whose elements are words in $a$, $b$, and
$(ab)$ which contain exactly $\ell$ letters $a$ and $r$ letters $b$. Letter $(ab)$ contributes one $a$ and one $b$. We say $X<Y$ if $X$ is obtained from $Y$ by adding parentheses or by replacing some number of $(ba)$'s by $(ab)$'s.
\end{definition}


For example, we have

\vspace{-0.55cm}

\begin{gather*}
aab\rightarrow a(ab)\leftarrow aba\rightarrow (ab)a\leftarrow baa\\
\text{Poset $X_{2,1}$}
\end{gather*}

\vspace{0.05cm}

\begin{prop}\label{p22}
For any $\ell,r\geq 0$, the poset $X_{\ell,r}$ is contractible.
\end{prop}

In order to prove this, we will embed the realization of $X_{\ell,r}$ in $\R^\ell$.

\begin{definition}\label{d23}
An {\it integer cube} of $\R^\ell$ is a cube of any dimension $s$, $0\leq s\leq \ell$, with  integer vertices and
whose edges are all of length one and parallel to one of the axes.
\end{definition}

The proof of the following is immediate.

\begin{lemma}\label{l24}
(i) An integer cube in $\R^\ell$ is determined by its center.
The dimension of an integer cube is the number of non-integer coordinates of the center.
The set (of centers) of integer cubes is the set ${\frac 12}\cdot\Z^\ell$ of points with half-integer coordinates.

(ii) $\R^\ell$ is a disjoint union of the interiors of integer cubes. Point $(x_{1},\ldots,x_{\ell})$ belongs to the
interior of the integer cube whose center has coordinates $(\frac{\lceil x_{1}\rceil+\lfloor x_{1}\rfloor}2, \ldots,
\frac{\lceil x_{\ell}\rceil+\lfloor x_{\ell}\rfloor}2)$.
\end{lemma}

\begin{figure}
\psfrag{x1}[0][0][1][0]{$x_1$}
\psfrag{x2}[0][0][1][0]{$x_2$}
\psfrag{X22}[0][0][1][0]{$|X_{2,2}|$}
\psfrag{X23}[0][0][1][0]{$|X_{2,3}|$}
\includegraphics[width=10cm]{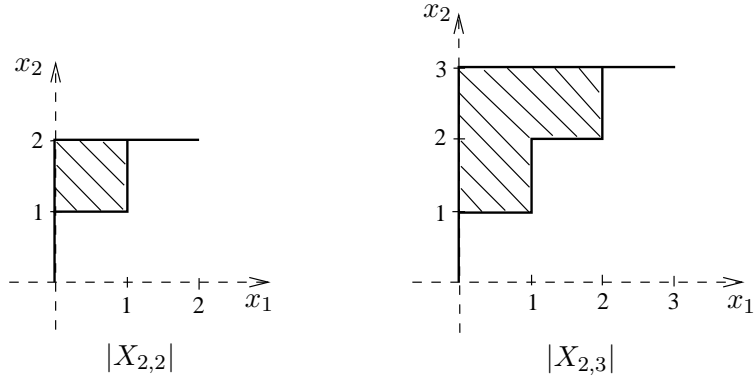}
\caption{Examples of $|X_{\ell,r}|$.}\label{Xlr}
\end{figure}

\begin{lemma}\label{l25}
The realization $|X_{\ell,r}|$ of $X_{\ell,r}$ is homeomorphic to the union of the integer cubes in
$\R^\ell=\{\,(x_{1},\ldots,x_{\ell})\,\}$ that are contained in the domain
$$
0\leq x_{1}\leq x_{2}\leq\ldots\leq x_{\ell}\leq r. \eqno(\numb)\label{eq21}
$$
\end{lemma}
Examples of $|X_{\ell,r}|$ are given in Figure \ref{Xlr}.

\begin{remark}\label{r26}
{\rm It follows from Lemma~\ref{l24}~(ii), that the subspace of $\R^\ell$ described above is defined by the
inequalities $\lceil x_{i-1}\rceil\leq\lfloor x_{i}\rfloor$, $i=1\ldots\ell+1$, where $x_{0}=0$ and $x_{\ell+1}=r$. }
\end{remark}

\begin{proof}  We define an embedding
$f\colon |X_{\ell,r}|\hookrightarrow\R^\ell$ on the elements of $X_{\ell,r}$. A simplex $X_{0}<X_{1}<\ldots<X_{k}$, $X_{i}\in X_{\ell,r}$ will then be mapped to the convex hull
of $f(X_{0}),\, f(X_{1}),\ldots,\, f(X_{k})$.

Let $X\in X_{\ell,r}$ be a word. We define $f(X)=(f_{1}(X),\ldots,f_{\ell}(X))$ as follows.
The coordinate $f_{i}(X)$ is set to be the number of elements $b$ before the $i$th $a$ in
$X$, but if this $a$ is parenthesized with a $b$, then  $\frac 12$ is added. For example,
$$
f(\,ab(ab)bbab\,)=(0,\frac 32, 4).
$$
The number of non-integer coordinates is exactly the number of letters $(ab)$ in $X$.

Note that the words in $X_{\ell,r}$ without parentheses are in one-to-one correspondence (via $f$) with the
integer points of the domain~\eqref{eq21}.  Similarly a half-integer point $N=(\frac{n_{1}}2,\frac {n_{2}}2,\ldots,\frac {n_{\ell}}2)$  is in the image of $f$ if and only
if the integer cube whose center is $N$ is contained in the domain~\eqref{eq21}.  Consider the
full subcategory $X_{\ell,r}\downarrow X$ of elements greater than or equal to some $X\in X_{\ell,r}$.  The realization
$|X_{\ell,r}\downarrow X|$ is homeomorphic to the barycentric subdivision of a cube (whose dimension is the number
of letters $(ab)$ in $X$). Space $|X_{\ell,r}\downarrow X|$ is mapped by $f$ to the integer cube with the center
$f(X)$.  Thus $f$ is an embedding, and the image is exactly the space described in the lemma.
\end{proof}

\begin{proof}[Proof of Proposition~\ref{p22}] We induct over $\ell$. The poset $X_{0,r}$ is
 contractible since it is a point. Consider the description of $|X_{\ell,r}|$ given in Remark~\ref{r26}.  The projection of $|X_{\ell,r}|$ to the first $\ell-1$ coordinates of $\R^{\ell}$ gives $|X_{\ell-1,r}|$.
The preimage of any point $(x_{1},\ldots,x_{\ell-1})\in f(|X_{\ell-1,r}|)$ is the segment
$[\lceil x_{\ell-1}\rceil,r]$ (in the degenerate case $\lceil x_{\ell-1}\rceil=r$, this segment is a point).
So $f(|X_{\ell,r}|)$ can be retracted to $f(|X_{\ell,r}|)\cap\{x_{\ell}=r\} \simeq |X_{\ell-1,r}|$,
which is contractible by the inductive hypothesis.
\end{proof}


\begin{definition}\label{d31}
{\rm Let $Y\in\T([n+1])$ and define $\Pi_n^{-1}(Y)$ to be the full subcategory of $\Phi(\nnn)$ with elements $\widehat Y$ satisfying
$\Pi_n(\widehat Y)=Y$.  }
\end{definition}

\begin{lemma}\label{l32}
Let $Y\in\T([n+1])$ have $\ell$ left-most nodes and $r$ right-most nodes (see Definition~\ref{d12}). Then $\Pi_n^{-1}(Y)$ is isomorphic to $X_{\ell,r}$.
\end{lemma}

\begin{proof} Suppose $\widehat Y\in\Pi_n^{-1}(Y)$.
We assign to $\widehat Y$ a word in letters $a$, $b$, $(ab)$ as follows:  If we travel along the trunk (see Definition~\ref{d19}) from the root to the distinguished leaf and meet a node that has only left-going branches, we write $a$.
If we meet a node that has only right-going branches, we write $b$. If this node has both left-going and right-going
branches, we write $(ab)$. Proceeding like this, we get a word in $X_{\ell,r}$.
(For example, the fan from Figure~\ref{fig0} produces $ba(ab)$.)
It is easy to see that such words are in one-to-one correspondence with the elements of $\Pi_n^{-1}(Y)$.
\end{proof}

The following is a consequence of Proposition~\ref{p22} and Lemma \ref{l32}.

\begin{cor}\label{c33}
$\Pi_n^{-1}(Y)$ is contractible for any $Y\in\Psi([n+1])$.
\end{cor}

\begin{definition}\label{d34}
Suppose given $Y\in\T([n+1])$ and $\widehat T\in\Phi(\nnn)$ such that $Y\geq\Pi_n(\widehat T)$. Define $\Pi_n^{-1}(Y|\geq \widehat T)$ as the
full subcategory of $\Phi(\nnn)$ whose elements $\widehat Y$ satisfy $\Pi_n(\widehat Y)=Y$ and $\widehat Y\geq \widehat T$.
\end{definition}

\begin{remark}\label{r35}
{\rm $\Pi_n^{-1}(Y)=\Pi_n^{-1}(Y| \geq *)$, where $*=\unitlength=0.26mm
\begin{picture}(30,20)
\put(10,15){\circle{5}} \put(10,0){\circle*{4}} \put(10,0){\line(0,1){12.5}} \put(8.,17){\line(-3,1){17}}
\put(12,17){\line(3,1){17}} \put(8.8,18.3){\line(-2,1){14}} \put(4,22){$\ldots$}
\end{picture}
$
 is the terminal (minimal) element  of $\Phi(\nnn)$. }
\end{remark}

\begin{prop}\label{p36}
Poset $\Pi_n^{-1}(Y| \geq \widehat T)$ is contractible for any $Y\in\T([n+1])$ and $\widehat T\in\Phi(\nnn)$ satisfying
$Y\geq\Pi_n(\widehat T)$.
\end{prop}

\begin{proof} We will prove that $\Pi_n^{-1}(Y| \geq \widehat T)$ is isomorphic
to $\prod_{i=0}^{s}X_{\ell_i,r_i}$ for some $s\geq 0$, $\ell_{i}, r_{i}\geq 0$, $i=0\ldots s$.  The result will follow from Proposition~\ref{p22}.

Label the vertices of the trunk of $\widehat T$ by $0,1,2,\ldots,s$ (from the root to the distinguished leaf), with $0$ corresponding to the root and $s$ to the last node before the distinguished leaf.  Define $\ell_{0}$
(resp. $r_{0}$) as the number of left-most (resp. right-most) vertices of $Y$ that are contracted to the root in $\Pi_n(\widehat T)$.
Analogously define $\ell_{i}, r_{i}$, $i=1,..., s$ as follows:  Denote by $I_{i}$, $i=1,..., s$, the $i$th vertex on the trunk (which is neither
the root nor the distinguished leaf). If $I_{i}$ does not have any left-going (resp. right-going) branches, then the set $\ell_{i}=0$ (resp. $r_{i}=0$).
Otherwise $I_{i}$ defines a left-most (resp. right-most) node $L_{i}$ (resp. $R_{i}$) of $\Pi_n(\widehat T)$. Then $\ell_{i}$ (resp. $r_{i}$)
is the number of left-most (resp. right-most) vertices of $Y$ contracted to $L_{i}$ (resp. $R_{i}$).




Poset $\Pi_n^{-1}(Y| \geq \widehat T)$ is a subposet of $\Pi_n^{-1}(Y)\cong X_{\ell,r}$, where $\ell=\sum_{i=0}^s\ell_{i},$ $r=\sum_{i=0}^s r_{i}$.  This subposet consists of those words in $X_{\ell,r}$ which can be broken up into $s+1$ words from $X_{\ell_{0},r_{0}}$, $X_{\ell_{1},r_{1}}$,..., $X_{\ell_{s},r_{s}}$.  But such a subposet is clearly the product $\prod_{i=0}^{s}X_{\ell_i,r_i}$.
\end{proof}


\begin{proof}[Proof of Theorem~\ref{t1}] The theorem follows from Corollary~\ref{c33}, Proposition~\ref{p36}, and a general proposition about maps of posets given below.
\end{proof}

\begin{prop}\label{p_contractible}
Let $\Pi\colon B\to A$ be a map of posets satisfying
\begin{itemize}
\item for any $Y\in A$, poset $\Pi^{-1}(Y)$ is contractible;\footnote{Posets $\Pi^{-1}(Y)$ and $\Pi^{-1}(Y|\geq \widehat T)$ below are defined similarly to Definitions~\ref{d31} and~\ref{d34}.}
\item for any $Y\in A$, and $\widehat T\in B$ such that $Y\geq \Pi(\widehat T)$, poset $\Pi^{-1}(Y|\geq \widehat T)$ is contractible;
\end{itemize}
then for any point of the realization $|A|$, its preimage under the induced map $|\Pi|\colon |B|\to|A|$ is contractible.
\end{prop}

\begin{proof} Let $z=z_{k}$ be a point in the open simplex $\Delta^k\subset |A|$ defined by a sequence
$$
Y_{0}<Y_{1}<\ldots<Y_{k},
\eqno(\numb)\label{eq41}
$$
where $Y_{i}\in  A$.
So $z_{k}$ is of the form
$$
z_{k}=\sum_{i=0}^{k}t_{i}Y_{i}; \quad \sum_{i=0}^kt_{i}=1; \quad t_{i}>0,\, i=0\ldots k.
$$
We will prove that $|\Pi|^{-1}(z_{k})$ is contractible by induction on $k$.

If $k=0$, then $|\Pi|^{-1}(z_{0})$ is the realization of $\Pi^{-1}(Y_{0})$, which is contractible by the first hypothesis.

Let $k\geq1$. Then $|\Pi|^{-1}(z_{k})$ has a natural prismatic decomposition.  Indeed, any $\widehat z_{k}\in|\Pi|^{-1}(z_{k})$ must be in the interior of some simplex
$$
(\widehat Y_{0}^0<\widehat Y_{0}^1<\ldots<\widehat Y_{0}^{m_{1}})<(\widehat Y_{1}^0<\ldots<\widehat Y_{1}^{m_{2}})<\ldots<
(\widehat Y_{k}^0<\ldots<\widehat Y_{k}^{m_{k}}),
\eqno(\numb)\label{eq42}
$$
$$
\widehat z_k = \sum_{{i=0\ldots k}\atop {j=0\ldots m_i}}t_i^j\widehat Y_i^j,
$$
where $m_{i}\geq 0$, $i=0,..., k$, $\Pi(\widehat Y_i^j)=Y_i$, and $t_{i}^j>0$ satisfy the condition
$\sum_{j=0}^{m_{i}}t_{i}^j=t_{i}$. (Parentheses are meant to simplify the reading of the expression.)
The sequence~\eqref{eq42} defines an open prism $\Delta^{m_{1}}\times\ldots\times\Delta^{m_{k}}$ in $|\Pi|^{-1}(z_{k})$. A prism is in the boundary of another if the sequence which defines the first is a subsequence of the sequence
which defines the second (however not any subsequence of~\eqref{eq42} defines a prism in
$|\Pi|^{-1}(z_{k})$; a subsequence that does has to contain at least one element $\widehat Y_{i}^j$ for each $i=0\ldots k$).

The preimages $|\Pi|^{-1}(z_{k})$ are naturally homeomorphic for all
$z_{k}\in\Delta^k\subset
|A|$ (recall that $\Delta^k$ is defined by the sequence~\eqref{eq41}). So
$|\Pi|^{-1}(\Delta^k)=\Delta^k\times|\Pi|^{-1}(z_{k})$.

Now let $z_{k-1}$ be any point in the open simplex $\Delta^{k-1}\subset |A|$ defined by the
sequence
$$
Y_{0}<Y_{1}<\ldots<Y_{k-1}
$$
and consider $|\Pi|^{-1}(z_{k-1})$. By induction hypothesis $|\Pi|^{-1}(z_{k-1})$ is
contractible. There is a natural map
$$
p\colon |\Pi|^{-1}(z_{k}) \longrightarrow |\Pi|^{-1}(z_{k-1})
$$
which is geometrically a boundary limit map. This is well-defined since $|\Pi|^{-1}(z_{k})$ is always the same space as $z_{k}\in\Delta^k$ tends to $z_{k-1}\in\Delta^{k-1}\subset \partial\Delta^k$. We will show that $p^{-1}$ of any point in $|\Pi|^{-1}(z_{k-1})$ is contractible.

In terms of prismatic decomposition, $p$ forgets the last factor in
$\Delta^{m_{1}}\times\ldots\times\Delta^{m_{k-1}}\times\Delta^{m_{k}}$ by mapping it to
$\Delta^{m_{1}}\times\ldots\times\Delta^{m_{k-1}}$ corresponding to the sequence
$$
(\widehat Y_{0}^0<\widehat Y_{0}^1<\ldots<\widehat Y_{0}^{m_{1}})<(\widehat Y_{1}^0<\ldots<\widehat Y_{1}^{m_{2}})<\ldots<
(\widehat Y_{k-1}^0<\ldots<\widehat Y_{k-1}^{m_{k-1}}).
\eqno(\numb)\label{eq44}
$$

Let $\widehat z_{k-1}$ be a point of the open prism $\Delta^{m_{1}}\times\ldots\times\Delta^{m_{k-1}}\subset|\Pi|^{-1}(z_{k-1})$.
It is easy to see that $p^{-1}(\widehat z_{k-1})$ is exactly the realization of the poset $\Pi^{-1}(Y_{k}| \geq\widehat Y_{k-1}^{m_{k}-1})$.
(Informally, we need to consider all the \lq\lq prolongations" of the sequence~\eqref{eq44} to the sequence~\eqref{eq42}.) But the last poset is always contractible (by the second hypothesis of the proposition). Thus $|\Pi|^{-1}(z_{k})$ is surjectively mapped by $p$ to $|\Pi|^{-1}(z_{k-1})$. Space $|\Pi|^{-1}(z_{k-1})$
is contractible, and the preimage of $p$ is contractible for any point of  $|\Pi|^{-1}(z_{k-1})$. So $|\Pi|^{-1}(z_{k})$ is contractible as well.\footnote{We use a well-known fact that a  map of compact $CW$-complexes with cell-like (specific type of contractible) preimages of any point is always a homotopy equivalence, see for example~\cite[Corollary~1.3]{Lacher} or~\cite[Corollary~2]{DydKo}.}
\end{proof}

\section{Associahedra and cyclohedra as compactifications of configuration spaces}\label{s5}

\subsection{General construction}\label{ss51}

Let $M$ be a smooth manifold. Denote by $C_n(M)$ the configuration space of $n$ distinct points in $M$, namely
$$
C_n(M)=\{(x_1,x_2,\ldots,x_n)\, | \, x_i\in M;\, x_i\neq x_j, 1\leq i\neq j\leq n \}.
$$
By $C_n[M]$ we will denote the Axelrod-Singer compactification of $C_n(M)$~\cite{AS} which is the differential-geometric analogue of the Fulton-McPherson algebraic-geometric compactification~\cite{FM}. In short the difference is that it can be obtained as a sequence of \lq\lq differential-geometric" spherical blow-ups instead of \lq\lq algebraic-geometric" blow-ups in projective spaces. We will use Kontsevich's version of this construction~\cite{Konts-OMDQ} that was fully developed by D.~Sinha in~\cite{Sinha-MCCS}. Assuming that $M$ is a submanifold of~$\R^m$, $C_n[M]$ can be defined as the closure of $C_n(M)$ in
$$
M^{\times n}\times (S^{m-1})^{\frac{n(n-1)}2}\times [0,+\infty]^{n(n-1)(n-2)}
\eqno(\numb)\label{eq51}
$$
under the inclusion
$$
\alpha_n=\iota\times (\pi_{ij}|_{C_n(M)})_{1\leq i<j\leq n}\times (s_{ijk}|_{C_n(M)})_{1\leq i\neq j\neq k \leq n},
\eqno(\numb)\label{eq52}
$$
where $\iota$ is the inclusion $\iota\colon C_n(M)\hookrightarrow M^{\times n}$; the maps $\pi_{ij}\colon C_n(\R^m)\to S^{m-1}$, $1\leq i\neq j\leq n$, are defined as $\frac{x_j-x_i}{||x_j-x_i||}$, and
$$
s_{ijk}\colon C_n(\R^m)\to [0,+\infty]
\eqno(\numb)\label{eq53}
$$
are $\frac{||x_i-x_j||}{||x_i-x_k||}$, $1\leq i\neq j\neq k\leq n$. The smooth structure on $[0,+\infty]$ is induced by the homeomorphism $[0,+\infty]\stackrel{\simeq}{\longrightarrow} [0,1]$ sending $t\mapsto \frac t{t+1}$.

The strata of $C_n[M]$ are in one-to-one correspondence with rooted trees having $n$ leaves labeled by $1,2,\ldots,n$. The only restriction on the trees is that all the non-leaf vertices except the root should be of valence $\geq 3$. Such a tree encodes how the points of configurations collide with each other. For example, a typical point ($=$ infinitesimal configuration) of the stratum encoded by the tree
$$
\includegraphics[width=5cm]{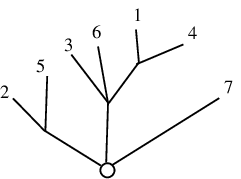}
$$
can be represented by the figure
$$
\includegraphics[width=12cm]{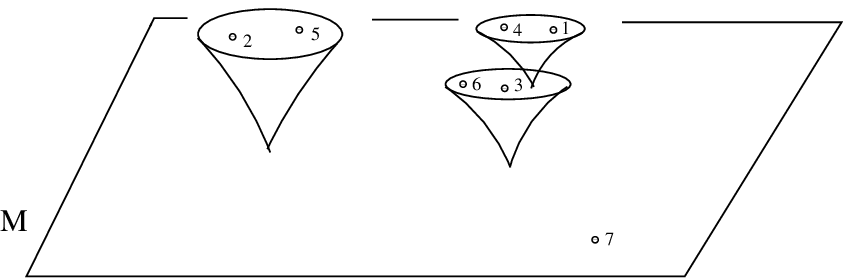}
$$
This figure tells us that the points 2 and 5 have the same image under the projection to $M$, and can be distinguished only due to the factor $S^{m-1}$ labeled by the map $\pi_{2,5}$. Similarly the points 1, 3, 4, 6 also all have the same image in $M$ (or in other words they collide), but the distance between~1 and~4 is infinitely small compared to the distances between~3, 6, and $1=4$. We will be assuming that the reader is familiar with this construction. As mentioned earlier, a very clear and detailed reference is the paper~\cite{Sinha-MCCS} by D.~Sinha.

One should mention that the last factor in~\eqref{eq51} is necessary to distinguish the limits of \lq\lq collinear" configurations, where the colliding points lie on the same line. In particular if $M$ is one-dimensional (the case we are interested in the most), all configurations are collinear.

\subsection{Associahedron}\label{ss52}
Consider the configuration space of $n+2$ points in $[0,1]$
$$
0=t_0<t_1<t_2<\ldots<t_n<t_{n+1}=1,
$$
where the first and the last points are fixed. Let $\Delta^n$ (open simplex) denote this space. The $n$-dimensional associahedron $Assoc_n$ can be defined as the closure of $\Delta_n$ under the inclusion
$$
\xymatrix{
\Delta^n\ar@{^{(}->}[r]^-{\alpha_n}& [0,1]^{{n+2}\choose 3} =: A_n,
}
\eqno(\numb)
\label{eq54}
$$
where $\alpha_n$ in coordinates is $s_{ijk}=(\alpha_n(t))_{ijk}=\frac{t_j-t_i}{t_k-t_i}$, $0\leq i<j<k\leq n+1$ (here $t$ denotes $(t_1,t_2,\ldots,t_n)$). Notice that
$$
s_{0,j,{n+1}}=\frac{t_j-t_0}{t_{n+1}-t_0}=t_j.
\eqno(\numb)\label{eq55}
$$
This explains why the first factor of~\eqref{eq51} is omitted in the definition of $A_n$. The second factor (which should be $(S^0)^{{n+2}\choose 2}$) is not necessary and omitted since we consider the compactification only of one connected component of the configuration space where the points on the line appear according to their linear order.

\begin{remark}\label{r51}
One should also mention that we use the functions $s_{ijk}$ only for $0\leq i<j<k\leq n+1$. We do so because the other $s_{ijk}$ can be smoothly expressed from them. Assuming $i<j<k$, one has:
$$
s_{ikj}=\frac 1{s_{ijk}}, \quad s_{jik}=\frac {s_{ijk}}{1-s_{ijk}}, \quad s_{jki}=\frac{1-s_{ijk}}{s_{ijk}}, \quad
s_{kij}=\frac 1{1-s_{ijk}}, \quad s_{kji}=1-s_{ijk},
\eqno(\numb)\label{eq55'}
$$
which are all smooth maps $[0,1]\to[0,+\infty]$.
\end{remark}

Abusing notation, we will use $(s_{ijk})_{0\leq i<j<k\leq n+1}$ as coordinates for $A_n$. It follows from~\eqref{eq55} that the coordinates
$
0<s_{0,1,n+1}<s_{0,2,n+1}<\ldots<s_{0,n,n+1}<1
$
can be used for the interior of $Assoc_n$.

We will be using the ambient space $A_n$ containing $Assoc_n$ when we need to define any map from or onto the associahedron. 

The strata of $Assoc_n=K_{n+2}$ are described by the poset $\Psi([n+1])$. For example, the stratum corresponding to the tree
$$
\includegraphics[width=5cm]{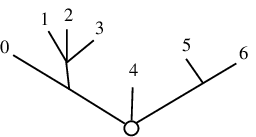}
$$
consists of infinitesimal configurations as follows:
$$
\includegraphics[width=10cm]{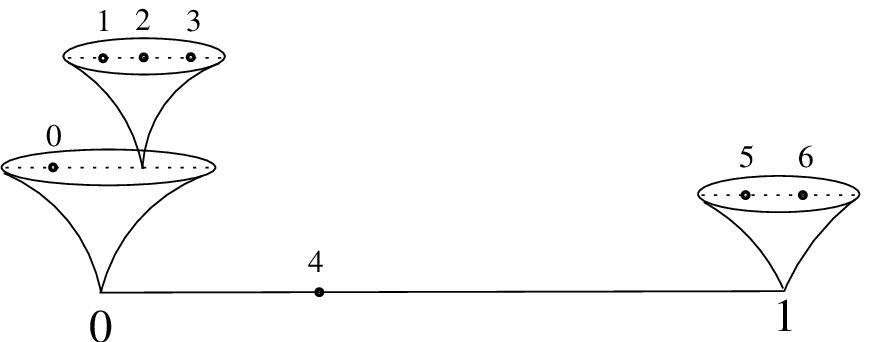}
$$
that are limits of configurations whose points 1, 2, 3 approach the left end of the interval and the point 5 approaches the right end of the interval. Moreover, the distance between 1 and 3 goes to zero infinitely faster than the distance between 0 and 1.

Let $T\in\Psi([n+1])$ be any tree. The face (closed stratum) $Assoc_T$ corresponding to $T$ is diffeomorphic to the product
$$
Assoc_T=\prod_{v\in V(T)}Assoc_{|v|-2},
\eqno(\numb)\label{eq56}
$$
where $V(T)$ is the set of internal vertices of $T$, and $|v|$ denotes the number of outgoing vertices of $v$. Let $e_1,\ldots,e_{|v|}$ be the outgoing edges of $v$ written in their natural linear order. We will use $0=t_{e_1}<t_{e_2}<\ldots<t_{e_{|v|-1}}<t_{e_{|v|}}=1$ as coordinates for the interior of the factor $Assoc_{|v|-2}$ in~\eqref{eq56}. Notice that for any edge $e$ of $T$ we assign a coordinate $t_e$, but for the minimal and maximal edges emanating from a vertex the corresponding coordinate is fixed to be $0$ and $1$ respectively.
The product~\eqref{eq56} is the closure of the image of the inclusion
$$
\xymatrix{
\prod_{v\in V(T)}\Delta^{|v|-2}\ar@{^{(}->}[r]^{\alpha_T}&\prod_{v\in V(T)} A_{|v|-2},
}
\eqno(\numb)\label{eq57}
$$
where $\alpha_T=\prod_{v\in V(T)}\alpha_{|v|-2}$. As coordinates on $A_{|v|-2}$ we will use $s_{e_i,e_j,e_k}$, where $e_i$, $e_j$, $e_k$ are edges outgoing from $v$ and satisfying $e_i<e_j<e_k$ (to recall the edges outgoing from a vertex have a natural linear order, see Section~\ref{ss21}).

To describe explicitly the inclusion $Assoc_T\hookrightarrow Assoc_n$ we need to define the inclusion $I_T$:

$$
\xymatrix{
\prod_{v\in V(T)}\Delta^{|v|-2}\ar@{^{(}->}[r]^{\alpha_T}&\prod_{v\in V(T)} A_{|v|-2}\ar@{^{(}->}[r]^-{I_T}&A_n
}.
\eqno(\numb)\label{eq58}
$$
The closure of the image of $I_T\circ\alpha_T$ is exactly the face $Assoc_T$.

\begin{definition}\label{d52}
The {\it nadir} $n(i,j)$ of leaves $i$ and $j$ of a tree $T$ is defined as the closest to the root vertex  on the path between $i$ and $j$.
\end{definition}

\begin{definition}\label{d53}
We will say that two leaves $i$ and $j$ are {\it closer} to each other than to a leaf $k$, if $n(i,j)\succ n(j,k)$ according to the natural partial order on the set $V(T)$ of vertices of $T$ (which simply means that $n(i,j)$ is above $n(j,k)$).
\end{definition}

In coordinates, the map $I_T$ is defined as follows:
$$
s_{ijk}=
\begin{cases}
0,&\text{if $i$ and $j$ are closer to each other than to $k$};\\
1,&\text{if $j$ and $k$ are closer to each other than to $i$};\\
s_{e(i),e(j),e(k)},&\text{if $n(i,j)=n(j,k)$}.
\end{cases}
\eqno(\numb)\label{eq58}
$$

\subsection{Cyclohedron}\label{s53}
The $n$-dimensional cyclohedron $Cycl_n$ was originally defined as the Axelrod-Singer compactification of the configuration space of $n+1$ points $(e^{2\pi it_0},e^{2\pi it_1},\ldots,e^{2\pi it_n})$, $0=t_0<t_1<\ldots<t_n<1$, on the circle $S^1$:
$$
\psfrag{R}[0][0][1][0]{$\mathrm{Re}\, z$}
\psfrag{I}[0][0][1][0]{$\mathrm{Im}\, z$}
\includegraphics[width=6cm]{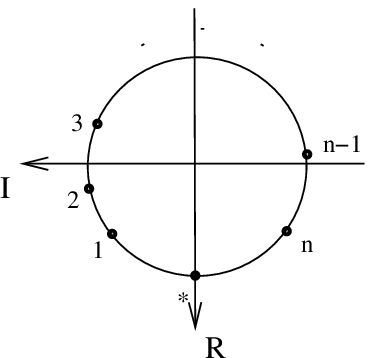}
$$
Notice that the first point $e^{t_0}=e^{t_*}=1$ (the distinguished point) is fixed. We used an unusual choice of coordinates in the above figure because we wanted the points to appear in  counter-clockwise order with the distinguished point at the bottom (like in Figure~\ref{homeofig}).

More precisely, the cyclohedron $Cycl_n$ is the closure of $\Delta^n$ under the inclusion
$$
\xymatrix{
\Delta^n\ar@{^{(}->}[r]^-{\gamma_n}&(S^1)^{n}\times (S^1)^{{n+1}\choose 2}\times [0,+\infty]^{(n+1)n(n-1)},
}
\eqno(\numb)\label{eq59}
$$
where $\gamma_n=(i_k)_{1\leq k\leq n}\times(\pi_{k\ell})_{0\leq k<\ell\leq n}\times (\tilde s_{k\ell m})_{0\leq k\neq \ell\neq m\leq n}$, with
$$
i_k(t)=e^{2\pi it_k},\qquad
\pi_{k\ell}(t)=e^{\pi i(t_k+t_\ell +\frac 12)},\text{ and} \quad
\tilde s_{k\ell m} (t)=\left| \frac {\sin \pi(t_\ell - t_k)}{\sin \pi(t_m-t_k)}\right|.
\label{eq510}
$$
Notice that $\pi_{kl}(t)$ is the unit vector
giving the direction from $i_k(t)$ to $i_\ell(t)$, and
$\tilde s_{k\ell m}(t)$ is a measure of the relative distance of
$i_k(t)$, $i_\ell(t)$, and $i_m(t)$.

Let $\bar\Delta^n$ denote the closed simplex $\{(t_1,t_2,\ldots,t_n)\,|\, 0\leq t_1\leq t_2\ldots\leq t_n\leq 1\}$, and let $\iota\colon\Delta^n\hookrightarrow\bar\Delta^n$ denote the inclusion of its interior.

\begin{lemma}\label{l51}
The cyclohedron can be defined as the closure of $\Delta^n$ under the inclusion 
$$
\xymatrix{
\Delta^n\ar@{^{(}->}[r]^-{\beta_n}&\bar\Delta^n\times [0,+\infty]^{(n+1)n(n-1)}=:B_n
}
\eqno(\numb)\label{eq511}
$$
where $\beta_n=\iota\times(\tilde s_{k\ell m})_{0\leq k\neq \ell\neq m\leq n}$.
\end{lemma}

\begin{proof}
Notice that the images of $\pi_{0k}=e^{\pi i(t_k+\frac 12)}$, $1\leq k\leq n$, form a configuration of $n$ points on the upper semicircle:
$$
\psfrag{R}[0][0][1][0]{$\Re z$}
\psfrag{I}[0][0][1][0]{$\Im z$}
\includegraphics[width=6cm]{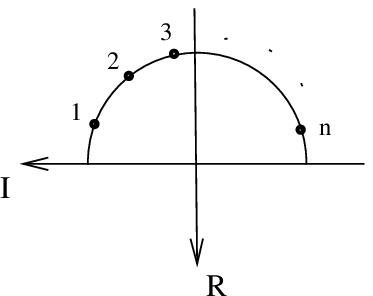}
$$
Therefore the closure of the images of $(\pi_{01},\pi_{02},\ldots,\pi_{0n})$ in $(S^1)^{\times n}$ is diffeomorphic to $\bar\Delta^n$. On the other hand, $\pi_{k\ell}$, $1\leq k<\ell\leq n$, and $i_k$, $1\leq k\leq n$, can be smoothly expressed in terms of $\pi_{0k}$, $1\leq k\leq n$. Indeed, $\pi_{k\ell}=e^{-\frac{\pi i}2}\pi_{0k}\cdot\pi_{0\ell}$, and $i_k=-(\pi_{0k})^2$. This means that we can get rid of the factors in~\eqref{eq59} corresponding to $i_k$, $1\leq k\leq n$, and $\pi_{k\ell}$, $1\leq k<\ell\leq n$.
\end{proof}

By abuse of notation we will use $(t_k)_{1\leq k\leq n}$, and $(\tilde s_{k\ell m})_{1\leq k\neq \ell\neq m\leq n}$ as coordinates on the space $B_n$.
Similarly to $A_n$ and $Assoc_n$, this space will be used to define any map from or onto $Cycl_n$.

The stratification of $Cycl_n$ is encoded by the poset $\Phi(\nnn)$. The following figure gives an example of a fan and an infinitesimal configuration from the corresponding stratum:
$$
\includegraphics[width=10cm]{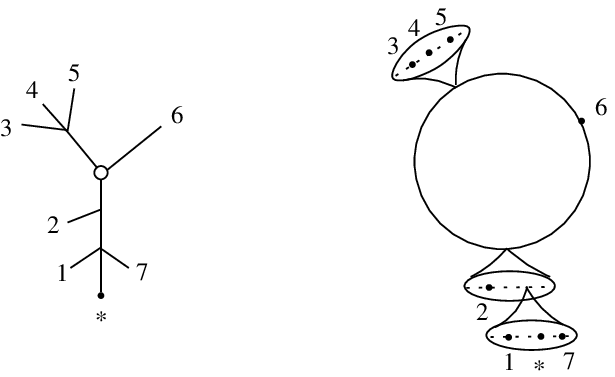}
$$
The right figure above represents an infinitesimal configuration which is the limit of configurations whose points 2 and 1 approach $*=0$ from the left and point 7 approaches it from the right, but the distance between 1 and 7 goes to zero infinitely faster than the distance between 2 and 1. The points 3, 4, and 5 also collide, the distance between 3 and 4 being comparable with the distance between 4 and 5.

The outgoing edges of any internal vertex (except the root) of a fan are linearly ordered. The outgoing edges of the root are cyclically ordered with one distinguished edge that points to the distinguished leaf of the fan.

For a fan $\widehat T\in\Phi(\nnn)$, denote by $Cycl_{\widehat T}$ the closed stratum of $Cycl_n$ corresponding to $\widehat T$. This face is diffeomorphic to the product
$$
Cycl_{\widehat T}=Cycl_{|root|-1}\times\prod_{v\in V(\widehat T)\setminus\{root\}}Assoc_{|v|-2},
\eqno(\numb)\label{eq512}
$$
where $root$ is the root, $V(\widehat T)$ is the set of non-leaf vertices of $\widehat T$ and $|\, .\, |$ denotes the number of outgoing edges of a vertex.

To define the inclusion $Cycl_{\widehat T}\hookrightarrow Cycl_n$ explicitly, one defines the inclusion $I_{\widehat T}$:
$$
\xymatrix{
\Delta^{|root|-1}\times\prod_{v\in V(\widehat T)\setminus\{root\}}\Delta^{|v|-2}\ar@{^{(}->}[r]^{\beta_{\widehat T}}& B_{|root|-1}\times\prod_{v\in V(\widehat T)\setminus\{root\}} A_{|v|-2}\ar@{^{(}->}[r]^-{I_{\widehat T}}& B_n,
}
\eqno(\numb)\label{eq513}
$$
where $\beta_{\widehat T}=\beta_{|root|-1}\times\prod_{v\in V(\widehat T)\setminus\{root\}}\alpha_{|v|-2}$. The closure of the image of $I_{\widehat T}\circ\beta_{\widehat T}$ is the face~\eqref{eq512}.

Below we describe $I_{\widehat T}$ in coordinates. Let $e_1$, $e_2,\ldots,$ $e_{|root|}$ be the edges outgoing from the root given in a cyclic order, with $e_1$ being the distinguished edge (pointing to the distinguished leaf). Then $0=t_{e_1}\leq t_{e_2}\leq t_{e_3}\leq\ldots\leq t_{e_{|root|}}\leq 1$, and $(\tilde s_{e_i,e_j,e_k})_{1\leq i\neq j\neq k\leq |root|}$ will be used as the coordinates of $B_{|root|-1}$. If $e_1$, $e_2$, $\ldots$, $e_{|v|}$ are the edges outgoing from a vertex $v\in V(\widehat T)\setminus\{r\}$ and that are given in the corresponding linear order, then $(s_{e_i,e_j,e_k})_{1\leq i< j< k\leq |v|}$ will be used as coordinates of $A_{|v|-2}$. The map $I_{\widehat T}$ in coordinates is given by~\eqref{eq513'} and~\eqref{eq512'}:
$$
t_i=t_{e(i)},\, 1\leq i\leq n,
\eqno(\numb)\label{eq513'}
$$
where $e(i)$ is the edge emanating from the $root$ that points to the leaf $i$. Recall Definitions~\ref{d52} and~\ref{d53}, one has
$$
\tilde s_{ijk}=
\begin{cases}
0,& \text{if $i$ and $j$ are closer to each other than to $k$};\\
1,& \text{if $j$ and $k$ are closer to each other than to $k$};\\
+\infty,& \text{if $i$ and $k$ are closer to each other than to $j$};\\
\tilde s_{e(i),e(j),e(k)},& \text{if $n(i,j)=n(j,k)=root$};\\
s_{e(i),e(j),e(k)},& \text{if $n(i,j)=n(j,k)\neq root$}.
\end{cases}
\eqno(\numb)\label{eq512'}
$$
In the above, $e(i)$, $e(j)$, and $e(k)$ are the edges outgoing from $n(i,j)=n(j,k)$ and pointing to the leaves $i$, $j$, and $k$ respectively. In the last case when $n(i,j)=n(j,k)$ is different from the root and one does not have  $e(i)<e(j)<e(k)$, then in addition one has to use~\eqref{eq55'} to define $s_{e(i),e(j),e(k)}$.

\section{Projection $\pi_n$}\label{s6}
In this section we define and study the projection $\pi_n\colon Cycl_n\to Assoc_n$ that was mentioned in the introduction. Recall that $Assoc_n$ is defined as the closure of $\Delta_n$ under the inclusion~$\alpha_n$:
$$
\xymatrix{
\Delta^n\ar@{^{(}->}[r]^-{\alpha_n}& [0,1]^{{n+2}\choose 3}\subset [0,+\infty]^{{n+2}\choose 3}=:A_n',
}
\eqno(\numb)\label{eq61}
$$
and $Cycl_n$ is the closure of the image of $\beta_n$:
$$
\xymatrix{
\Delta_n\ar@{^{(}->}[r]^-{\beta_n}& \bar\Delta^n\times [0,+\infty]^{(n+1)n(n-1)}= B_n.
}
\eqno(\numb)\label{eq62}
$$
$Assoc_n$ (respectively $Cycl_n$) is a smooth submanifold with corners of $A_n'$ (respectively $B_n$), see~\cite[Theorem~4.19]{Sinha-MCCS}. To define the projection $\pi_n\colon Cycl_n\to Assoc_n$, it is sufficient to provide a smooth map
$$
B_n \stackrel{p_n}{\longrightarrow} A_n'
\eqno(\numb)\label{eq63}
$$
satisfying $p_n\circ\alpha_n=\beta_n$. Let us define the $s_{ijk}$-component of the projection $p_n$.
First we notice that if $t_k-t_i>\frac 13$, then this component can be defined as
$$
s_{ijk}=f_{ijk}(t,\tilde s)=\frac{t_j-t_i}{t_k-t_i}.
\eqno(\numb)\label{eq64}
$$
In case $0\leq t_k-t_i<\frac 12$, it can be defined as
$$
s_{ijk}=g_{ijk}(t,\tilde s)=\tilde s_{ijk}\frac{\left(\sin\pi(t_k-t_i)\right)/(t_k-t_i)}{\left(\sin\pi(t_j-t_i)\right)/(t_j-t_i)}.
\eqno(\numb)\label{eq65}
$$
If $t_k-t_i=0$ or $t_j-t_i=0$, we take the limit. Since $0\leq t_j-t_i\leq t_k-t_i\leq\frac 12$, the right factor lies in $[\frac 2\pi,1]$. It is important that it is inside a strictly positive bounded interval otherwise the product with $\tilde s_{ijk}$ would not give a smooth function in $[0,+\infty]$. To define the $s_{ijk}$-component globally we use the partition of unity $\mu_1,\mu_2\colon\R\to[0,1]$ subordinate to the covering  $(-\infty,\frac 12)\cup (\frac 13,+\infty)$ of $\R$:
$$
s_{ijk}=\frac{\frac{\mu_1(t_k-t_i)f_{ijk}}{1+f_{ijk}}+\frac{\mu_2(t_k-t_i)g_{ijk}}{1+g_{ijk}}}
{1-\frac{\mu_1(t_k-t_i)f_{ijk}}{1+f_{ijk}}-\frac{\mu_2(t_k-t_i)g_{ijk}}{1+g_{ijk}}}.
\eqno(\numb)\label{eq66}
$$
The reason why the above formula is so complicated is that taking a linear combination in $[0,+\infty]$ would not give a smooth function, so one uses the diffeomorphism to [0,1], takes the linear (or more precisely affine) combination there, and then pulls the result back to $[0,+\infty]$.

There is some freedom in the definition of $p_n$. However the projection $\pi_n\colon Cycl_n\to Assoc_n$ is uniquely defined since it is determined by the homeomorphism of interiors (any continuous map is uniquely determined by its restriction to any dense subset). Later we will be assuming that $\mu_1$ and $\mu_2$ are always the same for all $s_{ijk}$-components, and for all $n$.

Once the projection $\pi_n\colon Cycl_n\longrightarrow Assoc_n$ is constructed we can formulate our second result.

\begin{theorem}\label{t2}
The preimage of any point of $Assoc_n$ under $\pi_n$ is contractible.
\end{theorem}

Theorem~\ref{t2} immediately follows from the following proposition.

\begin{proposition}\label{p61}
{\rm (i)} For any $\widehat T\in\Phi(\nnn)$, the image of the face $Cycl_{\widehat T}$ under the projection $\pi_n$ is the face $Assoc_T$, where $T=\Pi_n(\widehat T)$.

{\rm (ii)} The preimage of any point lying in the interior of $Assoc_T$, $T=\Pi_n(\widehat T)$, under the projection $\pi_{\widehat T}=\pi_n\bigr|_{Cycl_{\widehat T}}$, is homeomorphic to a $k$-dimensional cube, where $k$ is the number of vertices on the trunk of $\widehat T$ having both left and right outgoing edges.

{\rm (iii)} The preimage of any point lying in the interior of $Assoc_T$ under the projection $\pi_n$ is homemorphic to the realization of $\Pi_n^{-1}(T)$.
\end{proposition}

Recall from Lemma~\ref{l32} that $\Pi_n^{-1}(T)$ is isomorphic to $X_{\ell,r}$ for some $\ell$, and $r$. Thus (iii) together with Corollary~\ref{c33} imply Theorem~\ref{t2}. Notice also that (iii) gives an explicit geometric description of these preimages.

\begin{proof}
(i) follows from the commutative diagram~\eqref{eq68} that we consider below. A fan describes the way how the points collide in the infinitesimal configurations. Projection $\pi_n$ forgets how fast the points approach $*=e^{2\pi it_0}=1$ from the left compared to the points approaching $*$ from the right. The functor $\Pi_n$ does exactly the same thing.

(ii) Consider the restriction  $\pi_n\bigr|_{Cycl_{\widehat T}}=\pi_{\widehat T}$:
$$
\xymatrix{
Cycl_{\widehat T}\ar[rr]^(.62){\pi_{\widehat T}}\ar@{=}[d]&&Assoc_{T}\ar@{=}[d]\\
Cycl_{|root|-1}\times\prod\limits_{v\in V(\widehat T)\setminus\{root\}} Assoc_{|v|-2}\ar[rr]^-{\pi_{\widehat T}}&&\prod\limits_{v\in V(T)}Assoc_{|v|-2}.
}
\eqno(\numb)\label{eq67}
$$
To make this map explicit in coordinates, one has to consider the commutative diagram
$$
\xymatrix{
\Delta^{|root|-1}\times\prod_{v\in V(\widehat T)\setminus\{root\}}\Delta^{|v|-2}\ar[r]^-{\stackrel{\circ}{\pi}_{\widehat T}}
\ar@{^{(}->}[d]^{\beta_{\widehat T}}& \prod_{v\in V(T)}\Delta^{|v|-2}\ar@{^{(}->}[d]^{\alpha_{T}}\\
B_{|root|-1}\times\prod_{v\in V(\widehat T)\setminus\{root\}}A'_{|v|-2}\ar[r]^-{p_{\widehat T}}\ar@{^{(}->}[d]^{I_{\widehat T}}&
\prod_{v\in V(T)}A'_{|v|-2}\ar@{^{(}->}[d]^{I_T}\\
B_n\ar[r]^(.65){p_n}&A'_n
}
\eqno(\numb)\label{eq68}
$$
The maps $\alpha_T$, $I_T$, $\beta_{\widehat T}$, $I_{\widehat T}$ are defined in Section~\ref{s5}; the map $p_n$ is defined in the beginning of this section. The maps $\stackrel{\circ}\pi_{\widehat T}$, $p_{\widehat T}$ will be defined later in~\eqref{eq614}-\eqref{eq615}.
It is important to mention that the lower square of~\eqref{eq68} commutes because the functions $\mu_1$, $\mu_2$ in~\eqref{eq66} are the same
for all $i$, $j$, $k$, $n$, $0\leq i<j<k\leq n+2$. The above diagram is used to show that the projection $\pi_{\widehat
T}=\pi_n|_{Cycl_{\widehat T}}$ is a product of maps labeled by the non-leaf vertices of $\widehat T$:
$$
\pi_{\widehat T}=\pi_{|root|-1}\times\prod_{v\in V(\widehat T)\setminus\{root\}\setminus U(\widehat T)}{\mathrm id}_{Assoc_{|v|-2}}
\times\prod_{v\in U(\widehat T)}p_{r(v),\ell(v)},
\eqno(\numb)\label{eq69}
$$
where $U(\widehat T)$ is the set of those vertices on the trunk of $\widehat T$ that have both left and right outgoing edges; $r(v)$
(respectively $\ell(v)$) is the number of right (respectively left) outgoing edges of a vertex $v\in U(\widehat T)$; and the projection
$$
p_{r,\ell}\colon Assoc_{r+\ell-1}\longrightarrow Assoc_{r-1}\times Assoc_{\ell-1}
\eqno(\numb)\label{eq610}
$$
is specified below.

The idea of the factorization~\eqref{eq69} is that to any vertex of $\widehat T$ there correspond either one or two vertices of
$T=\Pi_n(\widehat T)$. To the $root$ of a fan, $\Pi_n$ assigns the root of $T$, which explains the presence of the factor $\pi_{|root|-1}$. To any vertex $v\in
V(\widehat T)\setminus\{root\}\setminus U(\widehat T)$ there corresponds exactly one vertex of $T$. Such a vertex $v$ produces the factor ${\mathrm id}_{Assoc_{|v|-2}}$ in~\eqref{eq69}. And finally, to any vertex $v\in U(\widehat T)$ there correspond 2 vertices of $T$:
$$
\psfrag{r}[0][0][1][0]{$r(v)$}
\psfrag{l}[0][0][1][0]{$\ell(v)$}
\psfrag{P}[0][0][1][0]{$\Pi_n$}
\includegraphics[width=10cm]{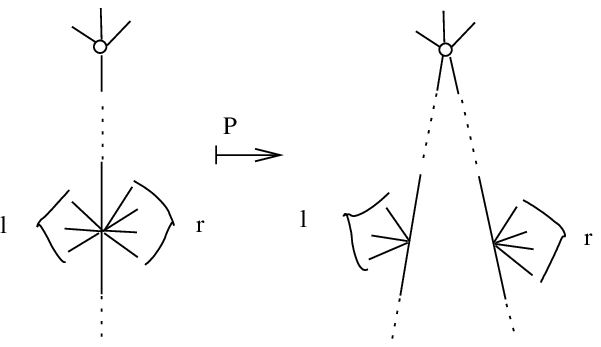}
$$
Such a vertex produces the factor
$$
p_{r(v),\ell(v)}\colon Assoc_{r(v)+\ell(v)-1}\longrightarrow Assoc_{r(v)-1}\times Assoc_{\ell(v)-1}.
\eqno(\numb)\label{eq611}
$$
The projection~\eqref{eq610} sends any point $0<t_1<t_2<\ldots<t_{r+\ell-1}<1$ in the interior of $Assoc_{r+\ell-1}$ to the pair
$$
\left(0<\frac{t_1}{t_{r}}<\frac{t_2}{t_{r}}<\ldots<\frac{t_{r-1}}{t_{r}}<1\, ; \,
0<\frac{t_{r+1}-t_{r}}{1-t_{r}}<\frac{t_{r+2}-t_{r}}{1-t_{r}}<
\ldots<\frac{t_{r+\ell-1}-t_{r}}{1-t_{r}}<1\right).
\eqno(\numb)\label{eq612}
$$
On the entire $Assoc_{r+\ell-1}$, the projection $p_{r,\ell}$ is defined as the restriction $\widehat p_{r,\ell}\bigr|_{Assoc_{r+\ell-1}}$, where $\widehat p_{r,\ell}$ is a map making the following diagram commute:
$$
\xymatrix{
\Delta^{r+\ell-1}\ar[r]^-{\stackrel{\circ}{p}_{r,\ell}}\ar@{^{(}->}[d]_{\alpha_{r+\ell-1}} & \Delta^{r-1}\times\Delta^{\ell-1}\ar@{^{(}->}[d]^{\alpha_{r-1}\times\alpha_{\ell-1}}\\
A_{r+\ell-1}\ar[r]^-{\widehat p_{r,\ell}}&A_{r-1}\times A_{\ell-1}}
\eqno(\numb)\label{eq613}
$$
The map $\stackrel{\circ}p_{r,\ell}$ is defined by~\eqref{eq612}. We use $(s^1_{ijk})_{0\leq i<j<k\leq r}$, and $(s^2_{ijk})_{0\leq i<j<k\leq \ell}$ as coordinates on $A_{r-1}$ and $A_{\ell-1}$ respectively. In coordinates, the map $\widehat p_{r,\ell}$ is given by
$$
\begin{matrix}
s^1_{ijk}=&s_{ijk},&0\leq i<j<k\leq r;\\
s^2_{ijk}=&s_{i+r,j+r,k+r},&0\leq i<j<k\leq \ell.
\end{matrix}
$$

\begin{lemma}\label{l62}
The preimage of any point in the interior of $Assoc_{r-1}\times Assoc_{\ell-1}$ under $p_{r,\ell}$ is a line segment.
\end{lemma}

\begin{proof}
Given any point $\tau\in[0,1]$, there is exactly one configuration (possibly infinitesimal if $\tau=0$ or 1) in $Assoc_{r+\ell-1}$ that projects to a given point in the interior of $Assoc_{r-1}\times Assoc_{\ell-1}$, and that has the propriety $s_{0,r,r+\ell}=\tau$ (see figure below).
$$
\psfrag{t0}[0][0][1][0]{$\tau=0$}
\psfrag{t1}[0][0][1][0]{$\tau=\frac 12$}
\psfrag{t2}[0][0][1][0]{$\tau=1$}
\includegraphics[width=15cm]{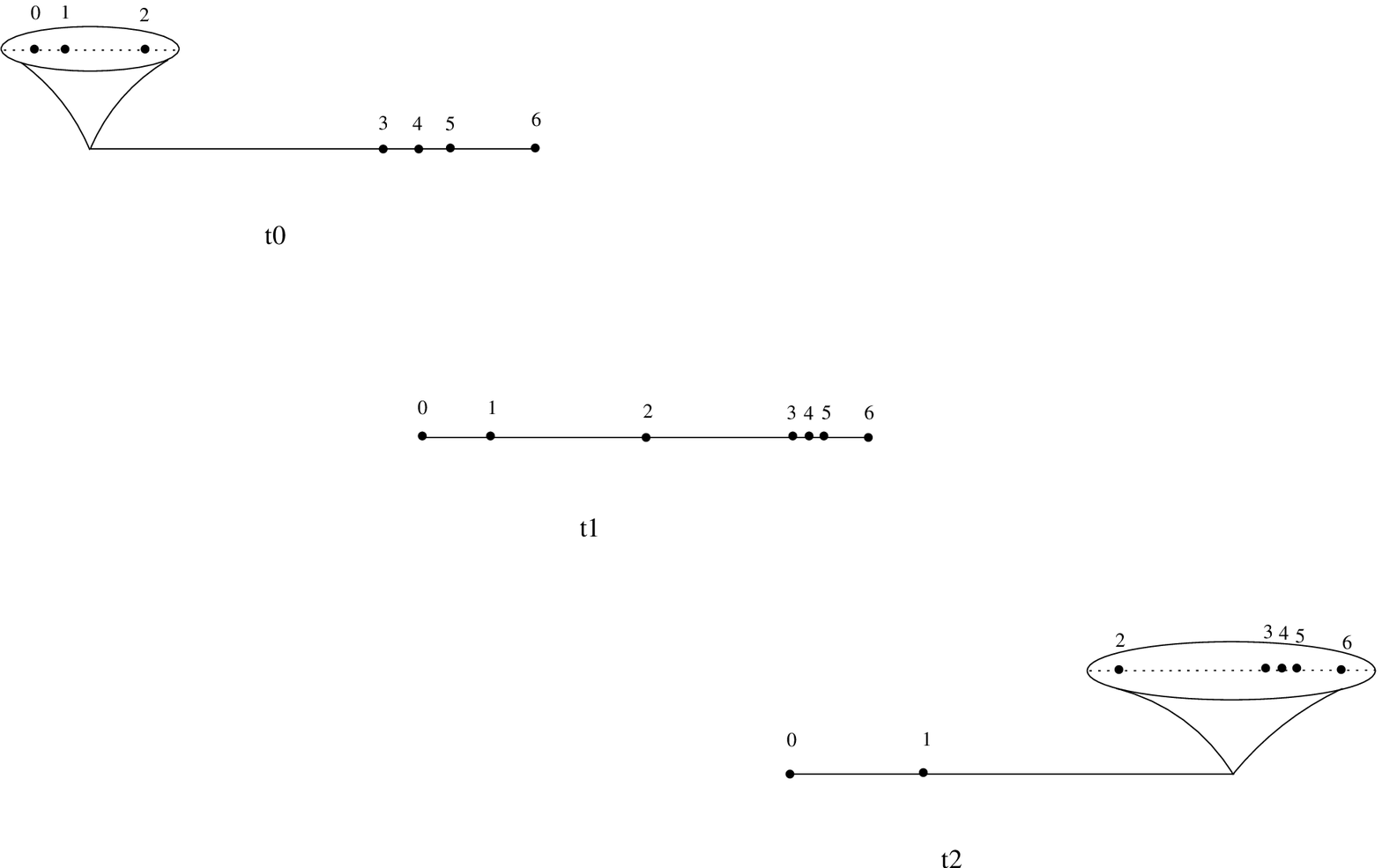}
$$
\end{proof}

Part (ii) of Proposition~\ref{p61} follows immediately from this lemma and from the factorization~\eqref{eq69}.\footnote{We also use the fact that $\pi_{|root|-1}$ is a homeomorphism of interiors.}

One should mention that one defines $\stackrel{\circ}\pi_{\widehat T}$, and $p_{\widehat T}$ in~\eqref{eq68} as follows:
$$
\stackrel\circ\pi_{\widehat T}=id_{\Delta^{|root|-1}}\times\prod_{v\in V(\widehat T)\setminus\{root\}\setminus U(\widehat T)}{\mathrm id}_{\Delta^{|v|-2}}
\times\prod_{v\in U(\widehat T)}\stackrel\circ p_{r(v),\ell(v)},
\eqno(\numb)\label{eq614}
$$
$$
p_{\widehat T}=p_{|root|-1}\times\prod_{v\in V(\widehat T)\setminus\{root\}\setminus U(\widehat T)}{\mathrm id}_{A'_{|v|-2}}
\times\prod_{v\in U(\widehat T)}{\widehat p}_{r(v),\ell(v)}.
\eqno(\numb)\label{eq615}
$$

(iii) Recall that $\Pi^{-1}(T)=X_{\ell,r}$, where $\ell$ (respectively $r$) is the number of left-most (respectively right-most) inner vertices of the tree $T$ (Lemma~\ref{l32}). The realization of $X_{\ell,r}$ has a natural cubical decomposition (see Lemma~\ref{l25}). Moreover, each cube corresponds to an element of $X_{\ell,r}$ which can be viewed as an element $\widehat T\in\Pi^{-1}(T)$. On the other hand, as it follows from (ii), the preimage of any point in the interior of $Assoc_T$ also has a natural cubical decomposition, the cubes being encoded by fans $\widehat T\in\Pi^{-1}(T)$. It follows that the two cubical complexes are isomorphic.
\end{proof}

\section{Permutohedron and leveled trees}\label{s7}

The $n$-dimensional permutohedron $Perm_n$ is defined as a convex hull in $\R^{n+1}$ of $(n+1)!$ points $(t_{\sigma(1)},t_{\sigma(2)},\ldots,
t_{\sigma(n+1)})$, $\sigma\in S_{n+1}$, where $t_1$, $t_2,\ldots,$~$t_{n+1}$ are distinct points of $\R$. The faces of $Perm_n$ are encoded by the ordered partitions of $\{1,2,\ldots,n+1\}$. It was shown in~\cite{Tonks-RAP} that this poset is isomorphic to the poset of planar leveled trees. Consider any planar leveled tree with $n+2$ leaves, and write  the numbers $1,2,\ldots,n+1$ between the leaves:

\vspace{0.5cm}

$$
\includegraphics[width=7cm]{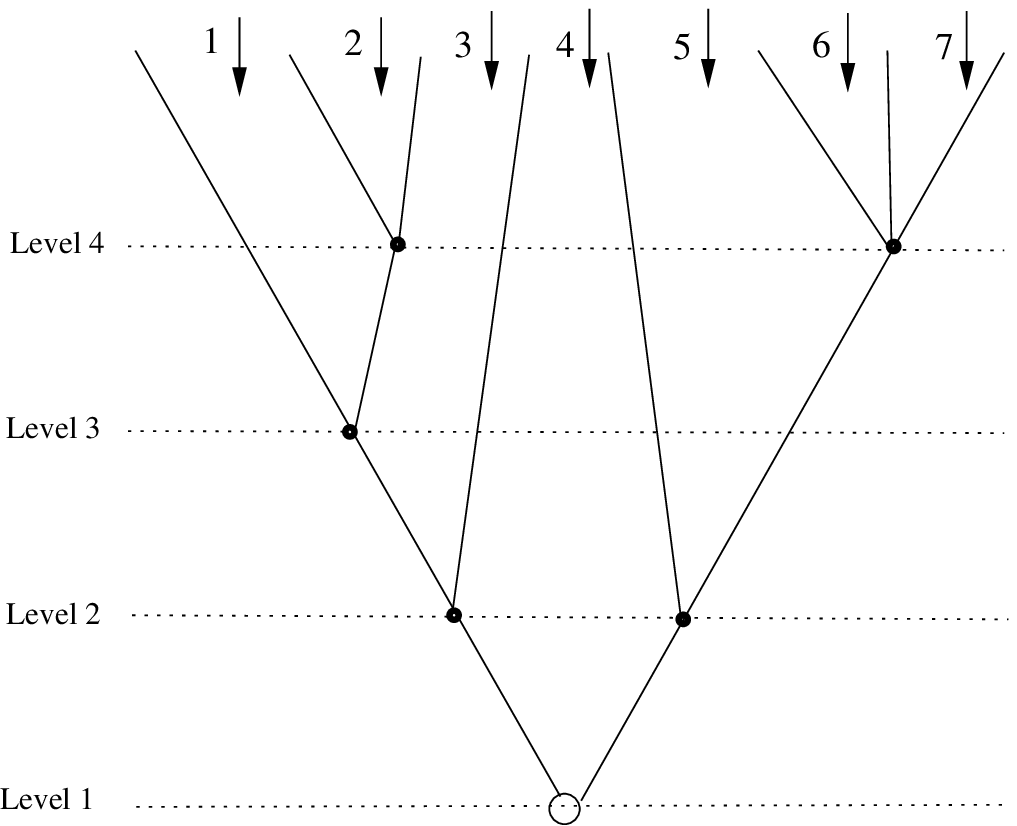}
$$

\vspace{1cm}

\noindent Then let the numbers fall down on the lowest possible vertex. The first set of the corresponding ordered partition is formed by the numbers that fell to the root (level~1), the second set by the numbers that fell to the second level, and so on. For example, for the leveled tree as above, the corresponding ordered partition is $\{4\}$, $\{3,5\}$, $\{1\}$, $\{2,6,7\}$.

In particular the leveled binary trees with one vertex on each level correspond by this construction to the permutations, or in other words encode the vertices of $Perm_n$.

To be precise, by a {\it leveled tree} in the sequel we will mean a rooted (not necessary planar) tree $T$ together with a surjective map $L$ from the set of its non-leaf nodes to some finite ordinal that respects the partial order $\prec$ induced by $T$.\footnote{For two nodes $a$ and $b$, one has $a\prec b$ if and only if the path from the root to $b$ passes through $a$.} We will say that a leveled tree $(T,L)>(T',L')$ if $(T',L')$ is obtained from $(T,L)$ by a contraction of levels. In particular $(T,L)>(T',L')$ implies $T\geq T'$.

Define a poset $\Psi^{level}([n+1])$ as a poset of leveled trees $(T,L)$ with $T\in\Psi([n+1])$. As mentioned earlier, this poset describes the face poset of $Perm_n$~\cite{Tonks-RAP}. The realization of $\Psi^{level}([n+1])$ is homeomorphic to the barycentric subdivision of $Perm_n$:
$$
|\Psi^{level}([n+1])|\simeq Perm_n.
\eqno(\numb)\label{eq71}
$$

Similarly we can define the poset $\Phi^{level}(\nnn)$ of leveled fans. Graphically the levels of a fan can be represented by concentric circles around the root:

\vspace{0.5cm}

$$
\includegraphics[width=6cm]{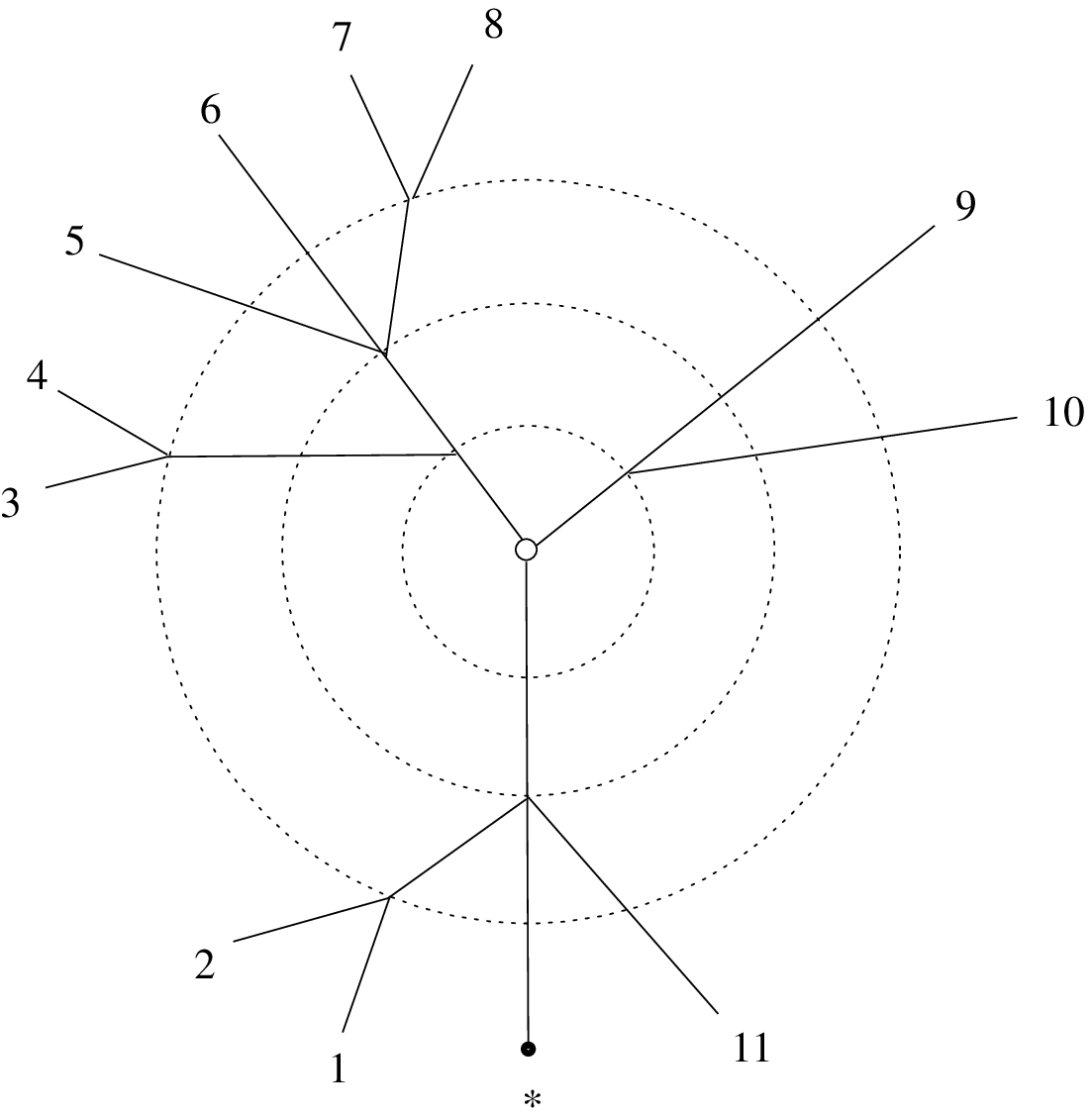}
$$

\vspace{.5cm}

\begin{lemma}\label{l71}
The posets  $\Phi^{level}(\nnn)$ and  $\Psi^{level}([n+1])$ are isomorphic.
\end{lemma}

\begin{proof}
The isomorphism $\Phi^{level}(\nnn)\stackrel{\simeq}{\longrightarrow}\Psi^{level}([n+1])$ is given by a similar surgery as the projection $\Pi_n\colon
\Phi(\nnn)\longrightarrow\Psi([n+1])$. For example,


$$
\includegraphics[width=12cm]{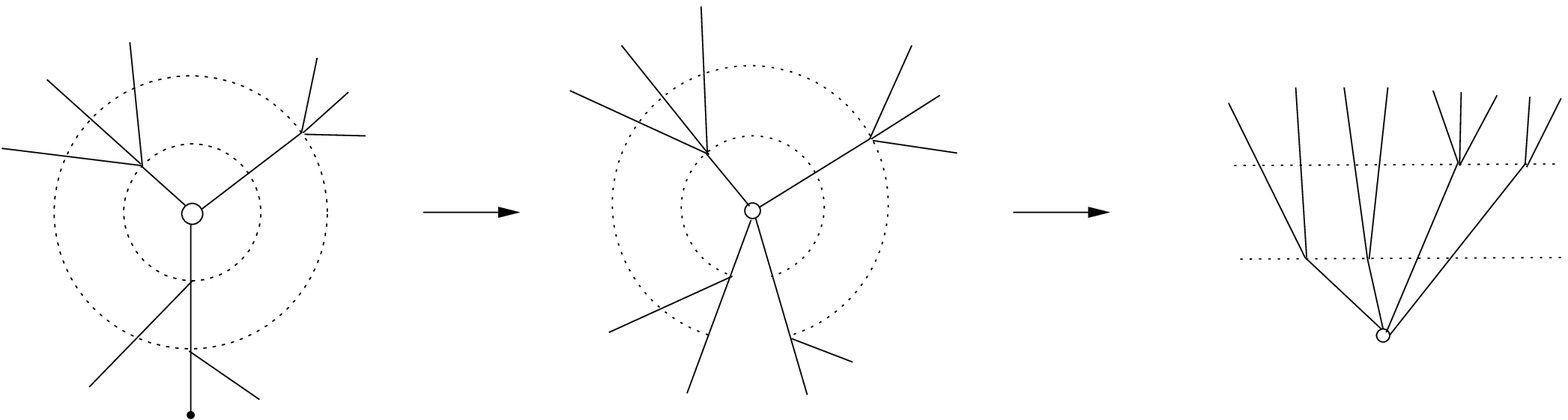}
$$
\end{proof}

Later on we will not distinguish between $\Phi^{level}(\nnn)$ and  $\Psi^{level}([n+1])$.

Forgetting the levels defines the projections
$$
\Phi^{level}(\nnn)\stackrel{\Pi'_n}\longrightarrow\Phi(\nnn) \hspace{.5cm} \text{and} \hspace{.5cm}
\Psi^{level}([n+1])\stackrel{\Pi''_n}\longrightarrow\Psi([n+1]).
\eqno(\numb)\label{eq72}
$$
The projection $\Pi''_n$ decomposes into a composition $\Pi''_n=\Pi_n\circ\Pi'_n.$ Passing to the realizations one gets maps
$$
Perm_n\stackrel{|\Pi'_n|}\longrightarrow Cycl_n \hspace{.5cm} \text{and} \hspace{.5cm}
Perm_n\stackrel{|\Pi''_n|}\longrightarrow Assoc_n.
\eqno(\numb)\label{eq74}
$$
Again, $|\Pi''_n|=|\Pi_n|\circ|\Pi'_n|$.

\begin{theorem}\label{t3}
For any point of $Cycl_n$ (respectively $Assoc_n$), its preimage under $|\Pi'_n|$ (respectively $|\Pi''_n|$) is  contractible.
\end{theorem}

The proof of this theorem is completely analogous to the proof of Theorem~\ref{t1}. It will follow from Proposition~\ref{p_contractible} and Lemmas~\ref{l72}-\ref{l73} below.

\begin{lemma}\label{l72}
For any rooted tree $T$, the poset $X_T$ of all possible levels on $T$ is contractible.
\end{lemma}

\begin{proof}
We induct on the number of non-leaf vertices of $T$. If $T$ has only one inner vertex (its root), $X_T$ is a point, and therefore contractible. Let $T$ have at least two inner nodes. Consider any of its non-leaf edges $e$. Define $T_e$ as the tree obtained from $T$ by contracting the edge $e$. One has a natural projection of posets
$$
X_T\stackrel{p_e}{\longrightarrow}X_{T_e}
$$
that \lq\lq forgets" the level of the upper vertex of the edge $e$. By induction hypothesis, we assume that $X_{T_e}$ is contractible. On the other hand, it is easy to see that the preimage of any point in $|X_{T_e}|$ under $|p_e|$ is homeomorphic either to a point or to a closed interval.\footnote{This is similar to the proof of Proposition~\ref{p22}, where one considers the projection $X_{\ell,r}\to X_{\ell-1,r}$. Notice by the way that the posets $X_{\ell,r}$ are examples of $X_T$ for some trees $T$.}
\end{proof}

\begin{lemma}\label{l73}
For any rooted tree $T$ and any leveled tree $(T',L')$ such that $T\geq T'$ the subposet of $X_T$ of levels $L$ on $T$ satisfying $(T,L)\geq (T',L')$ is contractible.
\end{lemma}

\begin{proof}
One can see that the poset in question is isomorphic to a product $\prod_iX_{T_i}$ taken over the set of levels of $(T',L')$, where each tree $T_i$ is obtained from $T$ by a contraction of some edges.
\end{proof}

\section{Leveled compactification of congiguration spaces}\label{s8}
There is another way to define $Cycl_n$ that makes the projection $\pi_n\colon Cycl_n\to Assoc_n$ more transparent.

\begin{proposition}\label{p81}
The $Cycl_n$ defined in Lemma~\ref{l51} is diffeomorphic as a manifold with corners to the closure of $\Delta_n$ under the inclusion
$$
\xymatrix{
\Delta^n\ar@{^{(}->}[r]^-{\gamma_n}&[0,1]^{{n+2}\choose 3}\times [0,+\infty]^{n\choose 2}=:C_n,
}
\eqno(\numb)\label{eq81}
$$
specified below.
\end{proposition}

In coordinates the map $\gamma_n$ has the same first ${n+2}\choose 3$ components as $\alpha_n$:
$$
s_{ijk}=\frac{t_j-t_i}{t_k-t_i}.
\eqno(\numb)\label{eq82}
$$
The last $n\choose 2$ components are
$$
r_{ij}=\frac{t_i-t_0}{t_{n+1}-t_j}=\frac{t_i}{1-t_j},\qquad 1\leq i<j\leq n.
\eqno(\numb)\label{eq83}
$$

\begin{proof}[Idea of the proof] One can proceed in the same way as in Section~\ref{s6} where the projection $\pi_n\colon Cycl_n\to Assoc_n$ was defined. Namely, one can define smooth maps $F_n$, $G_n$
$$
\xymatrix{
B_n\ar@/^/[r]^{F_n}&C_n,\ar@/^/[l]^{G_n}
}
\eqno(\numb)\label{eq84}
$$
such that $F_n\circ\beta_n=\gamma_n$ and $G_n\circ\gamma_n=\beta_n$ (where $B_n$, $\beta_n$ are defined in Section~\ref{s53}).
\end{proof}

The above model suggests that one should consider the closure of $\Delta_n$ under the inclusion
$$
\xymatrix{
\Delta^n\ar@{^{(}->}[r]^-{\delta_n}&[0,1]^{{n+2}\choose 3}\times [0,+\infty]^{{n+2}\choose 4}.
}
\eqno(\numb)\label{eq85}
$$
The first ${n+2}\choose 3$ components of $\delta_n$ are given by~\eqref{eq82}. The last ${n+2}\choose 4$ components are
$$
r_{ijk\ell}=\frac{t_j-t_i}{t_\ell-t_k},\qquad {0\leq i<j<k<\ell\leq n+1}.
\eqno(\numb)\label{eq86}
$$
It turns out that the closure of $\Delta_n$ under the inclusion $\delta_n$ is a manifold with corners homeomorphic to the $n$-dimensional permutohedron $Perm_n$.\footnote{This follows from Theorem~\ref{a82} analogously to the proof of~\cite[Theorem~4.19]{Sinha-MCCS}.}  This construction can be generalized to any manifold $M$. One defines $C_n[[M]]$, {\it leveled compactification} of $C_n(M)$, as the closure of $C_n(M)$ under the inclusion
$$
\xymatrix{
C_n(M)\ar@{^{(}->}[r]^-{\delta_n}&M^{\times n}\times (S^{m-1})^{n\choose 2}\times[0,+\infty]^{n(n-1)(n-2)}\times[0,+\infty]^{n(n-1)(n-2)(n-3)},
}
\eqno(\numb)\label{eq87}
$$
where $\delta_n=\alpha_n\times(r_{ijk\ell}\bigr|_{C_n(M)})_{1\leq i\neq j\neq k\neq \ell\leq n}$. The map $\alpha_n$ is defined by~\eqref{eq52} and
$$
\begin{array}{rccc}
r_{ijk\ell}\colon&C_n(\R^m)&\longrightarrow&[0,+\infty];\\
&(x_1,\ldots,x_n)&\mapsto&\frac{||x_i-x_j||}{||x_k-x_\ell||}
\end{array}.
\eqno(\numb)\label{eq87'}
$$

\begin{theorem-assertion}\label{a82}
The space $C_n[[M]]$ is a smooth submanifold with corners of the right side of~\eqref{eq87} whose strata are encoded by leveled trees with $n$ leaves labeled by
$1,2,\ldots,n$.
\end{theorem-assertion}

We give this result without proof. Even though technical it is pretty much straightforward and can be easily obtained by repeating the same arguments from~\cite{Sinha-MCCS}.

The levels appear due to the last factor of~\eqref{eq87}. Indeed, the functions $r_{ijk\ell}$ allow the comparison of the diameters of the distinct infinitesimal conglomerations  of points.\footnote{By \lq\lq conglomeration" we mean any colliding subset of points.} If one infinitesimal conglomeration has diameter 0 compared to another, the node of a tree corresponding to the first one appears on a level above the level containing the node corresponding to the other conglomeration. Notice that forgetting the last factor of~\eqref{eq87} induces a natural projection $C_n[[M]]\to C_n[M]$.


To convince the reader in the correctness of the above construction we should mention that there is an alternative method to obtain this leveled compactification which is analogous to the original Axelrod-Singer construction~\cite{Gaiffi}. Let $\calP$ be a partition of $\{1,2,\ldots,n\}$. We say that a point in $M^{\times n}=Maps(\{1,2,\ldots,n\},M)$ respects $\calP$ if it is constant inside each set of the partition $\calP$. Let $S_\calP(M^{\times n})$ denote the spherical blowup of $M^{\times n}$ along the strata of points respecting $\calP$. The leveled compactification $C_n[[M]]$ can be defined as the closure of $C_n(M)$ under the natural inclusion
$$
\xymatrix{
C_n(M)\ar@{^{(}->}[r]&\prod_\calP S_\calP(M^{\times n}),
}
\eqno(\numb)\label{eq89}
$$
where the product is taken over all partitions of $\{1,2,\ldots,n\}$. Using the general approach of Gaiffi~\cite{Gaiffi}, one can check that the strata of the obtained manifold with corners are indeed labeled by the leveled trees.
The difference from the Axelrod-Singer compactification $C_n[M]$ is that in the latter case the closure is taken in a subproduct of~\eqref{eq89} over those partitions of~$\{1,2,\ldots,n\}$ whose all elements except one are singletons. From this method we again recover the projection
$$
\xymatrix{
C_n[[M]]\ar[r]^{\pi_M}&C_n[M].
}
\eqno(\numb)\label{eq810}
$$

Denote by $C_T(M)$ the open stratum of $C_n[M]$ labeled by a tree $T$, and by $C_{(T,L)}(M)$ the open stratum of $C_n[[M]]$ labeled by a leveled tree $(T,L)$. We will denote their closures by $C_T[M]$ and $C_{(T,L)}[[M]]$.

\begin{theorem-assertion}\label{a83}
(i) The image of the stratum $C_{(T,L)}[[M]]$ under the projection $\pi_M$ is $C_T[M]$.

(ii) The preimage $\pi_M^{-1}(x)$ of any point $x\in C_T(M)$ inside $C_{(T,L)}[[M]]$ is homeomorphic to a product of permutohedra:
$$
\pi_M^{-1}(x)\cap C_{(T,L)}[[M]]\simeq\prod_i Perm_{k_i-1},
\eqno(\numb)\label{eq811}
$$
where the product is taken over the levels of $(T,L)$, and $k_i$ is the number of vertices of $(T,L)$ on the $i$-th level.

(iii) The preimage $\pi_M^{-1}(x)$ of any point $x\in C_T(M)$ is homeomorphic to the realization of the poset $X_T$ (see Lemma~\ref{l72}):
$$
\pi_M^{-1}(x)\simeq |X_T|.
\eqno(\numb)\label{eq811'}
$$
\end{theorem-assertion}

\begin{proof}[Sketch of the proof]
(i) follows from the fact that the image of $C_{(T,L)}(M)$ under $\pi_M$ is $C_T(M)$, which is straightforward from the geometric description of strata in terms of the collisions of points.

(ii) One can explicitly in coordinates define the commutative diagram~\eqref{eq812} below that defines the inclusion $\prod_i Perm_{k_i-1}\hookrightarrow C_n[[M]]$ lying inside  $C_{(T,L)}[[M]]$, and such that the interior of the product $\prod_i Perm_{k_i-1}$ maps exactly  to $\pi_M^{-1}(x)\cap{C}_{(T,L)}(M)$.
$$
\xymatrix{
\prod_i\Delta^{k_i-1}\ar@{^{(}->}[rr]^-{\prod_i\delta_{k_i-1}}\ar@{^{(}->}[d]_{i_1}&&\prod_i\left([0,1]^{{k_i+2}\choose 3}\times[0,+\infty]^{{k_i+2}\choose 4}\right)\ar[d]^{i_2}\\
C_{(T,L)}(M)\ar[d]_{\pi_M\bigr|_{C_{(T,L)}}}\ar@{^{(}->}[rr]^-{i_3}&&M^{\times n}\times(S^{m-1})^{n\choose 2}
\times [0,+\infty]^{n(n-1)(n-2)}\times [0,+\infty]^{n(n-1)(n-2)(n-3)}\ar[d]^p\\
C_{T}(M)\ar@{^{(}->}[rr]^-{i_4}&&M^{\times n}\times(S^{m-1})^{n\choose 2}
\times [0,+\infty]^{n(n-1)(n-2)}
}
\eqno(\numb)\label{eq812}
$$

In the above diagram, the product $\prod_i\Delta^{k_i-1}$  describes the interior of $\prod_i Perm_{k_i-1}$. The top right space in the diagram is the ambient space containing  $\prod_i Perm_{k_i-1}$; the one below it is the ambient space containing $C_n[[M]]$; and the lower right is the ambient space containing $C_n[M]$.

To get an idea about the map $i_2$, note that all its components except the last (corresponding to the factor $[0,+\infty]^{n(n-1)(n-2)(n-3)}$) are constant. Equivalently, $p\circ i_2$ is a constant map. For the inclusion $i_1$, consider as an example the leveled tree


$$
\includegraphics[width=4cm]{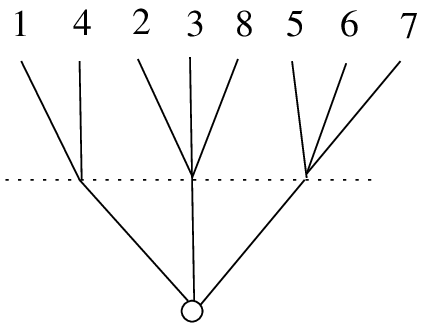}
$$

The product $\prod_i\Delta^{k_i-1}$ in this case has only one factor $\Delta^2$. The map $i_1$  assigns an infinitesimal leveled configuration in such a way that the ratio of distances
$$
d(x_1,x_4):d(x_2,x_3):d(x_5,x_6)
$$
is proportional to
$$
(t_1-t_0):(t_2-t_1):(t_3-t_2)=t_1:(t_2-t_1):(1-t_2).
$$

For (iii), first we notice that the realization $|X_T|$ is homeomorphic to the barycentric subdivision of a complex whose cells are labeled by the elements of $X_T$ and are products of permutohedra. Indeed, one can see that for any levelization $L$ of $T$, the subposet $X_T\downarrow L$ of levels that can be contracted to $L$ is isomorphic to $\prod_i\Psi^{level}([k_i-1])$, where $k_i$ is the number of vertices at the $i$-th level (to see this one should use the initial description of the poset of faces of a permutohedron in terms of ordered partitions; see Section~\ref{s7}). This means that the cell $|X_T\downarrow L|$\footnote{See Section~\ref{ss81} for the definition of $\calF\downarrow d$.} of $|X_T|$ is homeomorphic to $\prod_i Perm_{k_i-1}$. It follows from~(ii) that $\pi_M^{-1}(x)$ has exactly the same cell structure.
\end{proof}

Now we return to the permutohedron. There are two ways to define it as a compactification of a configuration space. One way, described at the beginning of the section, is the leveled compactification of the configuration space of $n+2$ linearly ordered points on $I=[0,1]$ with the first and the last points fixed at~$t=0$ and $t=1$ respectively; $Perm_n=C_n[[I,\partial]]$. The other way is to define $Perm_n$ as the leveled compactification of the configuration space of $n+1$ cyclically ordered points on a circle with one being fixed; $Perm_n=C_n[[S^1,*]]$. Notice that in the first case, the face poset is naturally described by $\Psi^{level}([n+1])$. In the second case the face poset is $\Phi^{level}(\nnn)$.  Similarly to the proof of Proposition~\ref{p81}, one can show that these two constructions produce diffeomorphic manifolds with corners, which explains in particular why the above two posets are isomorphic. The projections from the leveled compactifications to the usual ones produce maps
$$
Perm_n\stackrel{\pi'_n}{\longrightarrow}Cycl_n\hspace{.5cm} \text{and} \hspace{.5cm}
Perm_n\stackrel{\pi''_n}{\longrightarrow}Assoc_n.
\eqno(\numb)\label{eq813}
$$
One can easily see that $\pi''_n=\pi_n\circ\pi'_n$, since the maps coincide on the interior of $Perm_n$.

\begin{theorem}\label{a4}
The preimage of any point of $Cycl_n$ (respectively $Assoc_n$) under $\pi'_n$ (respectively $\pi''_n$) is contractible.
\end{theorem}

This theorem follows from Theorem~\ref{a84} and Lemma~\ref{l72}. The proposition below is immediate from Theorem~\ref{a83}~(iii).

\begin{proposition}\label{a84}
The preimage of any point inside the face $Cycl_{\widehat T}$, $\widehat T\in\Phi(\nnn)$, {\rm (}respectively $Assoc_T$, $T\in\Psi([n+1])${\rm )} under $\pi'_n$ (respectively $\pi''_n=\pi_n\circ\pi'_n$) is homeomorphic to $|X_{\widehat T}|=|(\Pi'_n)^{-1}(\widehat T)|$ {\rm (}respectively $|X_T|=|(\Pi''_n)^{-1}(T)|${\rm )}.
\end{proposition}

\section{Applications}\label{s_application}

\subsection{Cofinality}\label{ss81}

Let $\calF\colon\calC\to\calD$ be a functor between two small categories. Recall that, for any object
$d\in Ob(\calD)$, $\calF\downarrow d$ is defined as a category whose objects
are pairs $(c,f)$, where $c\in\calC$ and $f\in Mor_\calD(\calF(c),d)$. Morphisms are defined as
$$
Mor_{\calF\downarrow d}(\, (c_1,f_1);\, (c_2,f_2)\,) =
\{\, f\in Mor_\calC (c_1,c_2)\, |\, f_2\circ\calF (f)=f_1\,\}.
$$

Functor $\calF\colon\calC\to\calD$ is said to be {\it left cofinal} if, for any object
$d\in Ob(\calD)$, the realization of $\calF\downarrow d$ is contractible~\cite[Ch.~XI, page 316]{BK}.  This notion is important in homotopy theory since left cofinal functors preserve homotopy limits.  More precisely, for $X\colon \calD\to \Top$ a functor
from $\calD$ to the category of topological spaces and $\calF\colon\calC\to\calD$ left cofinal, we have
$$
\holim\limits_\calC X\circ\calF \simeq \holim\limits_\calD X.
$$


The following theorem was the main motivation for this paper. We discuss its applications
in  next subsection.

\begin{theorem}\label{t5}
The functors $\Pi_n$, $\Pi'_n$, and $\Pi''_n$ are left cofinal for any $n\geq 0$.
\end{theorem}

\begin{proof}
This result follows immediately from Theorems~\ref{t1} and~\ref{t3}. Consider, for example, the case of~$\Pi_n$.
We have to prove that for any $T\in\Psi([n+1])$, the realization of the category $\Pi_n\downarrow T$ is contractible. This realization is the preimage under $|\Pi_n|$ of the face of $Assoc_n$ encoded by $T$. We have that any face of $Assoc_n$, being a convex polytope,
is contractible. By Theorem~1 the preimage of any its point is also contractible.  Therefore using already mentioned~\cite[Corollary~1.3]{Lacher} or~\cite[Corollary~2]{DydKo} we conclude that the preimage of any face is also contractible.
\end{proof}

\begin{remark}
{\rm Projections $\Pi_n$, $\Pi'_n$, $\Pi''_n$ are clearly right cofinal. They preserve homotopy colimits
simply because they send the terminal object to the terminal object.}
\end{remark}

\subsection{Goodwillie-Weiss tower for spaces of knots}\label{ss43}
The space of long knots $Emb_d$, $d\geq 3$, is the space of smooth embeddings $f\colon
\R\hookrightarrow\R^d$ that coincide with a fixed linear embedding $t\mapsto (t,0,\ldots,0)$ outside
a compact subset of $\R$.
T.~Goodwillie and M.~Weiss defined a tower of spaces
$$
P_0Emb_d\leftarrow P_1Emb_d\leftarrow\ldots \leftarrow P_nEmb_d\leftarrow\ldots
\eqno(\numb)\label{eq52}
$$
converging to $Emb_d$ for $d\geq 4$ \cite{GW}.
Each space $P_nEmb_d$ is defined as a homotopy limit over a subcubical diagram, i.e.
over the category of faces of the $n$-simplex.

D.~Sinha gave several different models for $P_nEmb_d$ \cite{Sinha-TKS,Sinha-OKS}. In one of these models, the homotopy limit is taken over the category $\Psi([n+1])$ of faces of
$Assoc_n$. We will call this model {\it associahedral}.
On the other hand, the authors defined a {\it cyclohedral} model for the tower~\eqref{eq52} in~\cite{LTV} using the construction of a fanic diagram assigned to a morphism of operads.
%
The advantage of the last model is that the underlying fanic diagram is $\Q$-formal, which allows one to determine the rational homotopy type of $P_nEmb_d$ and $Emb_d$ for $d\geq 4$ \cite{LTV}.
The proof of this result does not directly use  Theorem~\ref{t2}. In that proof, it was enough that the composite map from $\Phi(\nnn)$ to the poset of faces of the $n$-simplex is left cofinal. However, Theorem~\ref{t2} can be used to give a more direct geometric relation between the
cyclohedral model from~\cite{LTV} and the associahedral model from \cite{Sinha-TKS}.  This is a part of the program aimed at relating Bott-Taubes integrals~\cite{BT} and their generalizations~\cite{CCRL1,CCRL2} to the Goodwillie-Weiss calculus of knots.  This program was partially realized by the third author in~\cite{Volic}
on the level of finite-type invariants of knots in~$\R^3$.

One can also define a {\it permutohedral} model for $P_nEmb_d$. The entries of the corresponding diagram are leveled compactifications of configuration spaces of points in $\R^d$ (or more generally in a manifold $M$ if we consider spaces of knots in $M$). This model completes the picture, though it does not seem to be of any immediate use. Theorem~\ref{t5} for the functors~$\Pi'_n$,~$\Pi''_n$ together with Theorem~\ref{a83} imply that the permutohedral model is equivalent to the cyclohedral and associahedral models constructed earlier.

\subsection*{Acknowledgements}
The authors are grateful to G.~Gaiffi, M.~Markl, D.~Sinha, and A.~Tonks for communications.

\end{document}